\newif\ifpreprint

\preprinttrue  

\ifpreprint

\documentclass[12pt]{article}

\usepackage{graphicx}

\usepackage{abstract} 

\usepackage{titlesec}
\def\sectionfont{\sffamily\Large\bfseries\boldmath}
\def\subsectionfont{\sffamily\large\bfseries\boldmath}
\def\paragraphfont{\sffamily\normalsize\bfseries\boldmath}
\titleformat*{\section}{\sectionfont}
\titleformat*{\subsection}{\subsectionfont}
\titleformat*{\subsubsection}{\paragraphfont}
\titleformat*{\paragraph}{\paragraphfont}
\titleformat*{\subparagraph}{\paragraphfont}
\usepackage[small,labelfont={bf,sf}]{caption}  

\usepackage{amsmath,amssymb,amsthm,enumitem}

\else

\documentclass[ijoc,nonblindrev]{informs3} 

\OneAndAHalfSpacedXII 

\usepackage{natbib}
 \bibpunct[, ]{(}{)}{,}{a}{}{,}%
 \def\bibfont{\small}%

\renewcommand{\cite}{\citep}

\TheoremsNumberedThrough     

\EquationsNumberedThrough    

\fi

\usepackage{mathtools}
\usepackage{multirow}
\usepackage{algorithm}
\usepackage[noend]{algpseudocode}
\algnewcommand{\IfThen}[2]{
  \State \algorithmicif\ #1\ \algorithmicthen\ #2}

\usepackage[
tight-spacing=true,
round-mode=places,
table-number-alignment=left,
table-figures-integer=1,
round-precision=2,
]{siunitx}
\usepackage{csvsimple}	
\usepackage{makecell, booktabs}   
\usepackage{adjustbox}  
\usepackage{adjustbox} 
\newlength{\tablelength}
\ifpreprint
\else
\fi

\usepackage{tikz}
\usetikzlibrary{calc}
\usetikzlibrary{positioning}
\definecolor{main}{RGB}{163, 31, 52}
\definecolor{secondary}{RGB}{4, 30, 62}
\definecolor{splitColor}{RGB}{202, 200, 200}
\definecolor{classColor1}{RGB}{163, 31, 52}
\definecolor{classColor2}{RGB}{4, 30, 65}
\definecolor{classColor3}{RGB}{0, 105, 143}
\definecolor{classColor4}{RGB}{4, 30, 65}
\tikzset{%
  every neuron/.style={
    circle,
    draw,
    minimum size=.5cm
  },
  neuron missing/.style={
    draw=none,
    scale=4,
    text height=0.333cm,
    execute at begin node=\color{black}$\vdots$
  },
}
\tikzset{class1/.style={draw,circle,text=white, fill=classColor1!80}}
\tikzset{class2/.style={draw,circle,text=white, fill=classColor2!80}}
\tikzset{class3/.style={draw,circle,text=white, fill=classColor3!80}}
\tikzset{class4/.style={draw,circle,text=white, fill=classColor4!80}}
\tikzset{split/.style={draw,rectangle,fill=splitColor!20}}
\tikzset{level 1/.style={sibling distance=12em},
         level 2/.style={sibling distance=6em},
         level 3/.style={sibling distance=3em}
        }

\usetikzlibrary{arrows}
\tikzstyle{block}=[draw, align=center, inner sep=.5em, fill=main, text=white, minimum size=2em]
\tikzstyle{init} = [pin edge={to-,thin,black}]

\usepackage{hyperref}

\usepackage{pgfplots, pgfplotstable}
\usepgfplotslibrary{statistics}
\usepgfplotslibrary{groupplots}
\usepgfplotslibrary{dateplot}
\usepgfplotslibrary{colorbrewer}

\usepackage[xindy, acronyms]{glossaries}   
\makeglossaries

\definecolor{borderBoxTodo}{rgb}{0,0.309,0.6}
\definecolor{backgroundBoxTodo}{rgb}{0.701,0.854,1}
\definecolor{lineColorTodo}{rgb}{0.701,0.854,1}
\RequirePackage[colorinlistoftodos, backgroundcolor=backgroundBoxTodo, bordercolor = borderBoxTodo, linecolor = lineColorTodo]{todonotes}

\newcommand{\reviewChanges}[1]{{#1}}

\ifpreprint

\newtheorem{theorem}{Theorem}[section]  

\fi

\newcommand{\eg}{{\it e.g.}}
\newcommand{\ie}{{\it i.e.}}

\newcommand{\ones}{\mathbf 1}
\newcommand{\reals}{{\mbox{\bf R}}}
\newcommand{\integers}{{\mbox{\bf Z}}}
\newcommand{\symm}{{\mbox{\bf S}}}  
\newcommand{\NPhard}{\mbox{$\mathcal{NP}$-hard}}  

\newcommand{\prob}{{\mathbf P}}

\newcommand{\identity}{I}
\newcommand{\tpose}{T}

\newcommand{\Expect}{\mathop{\bf E{}}}

\newcommand{\card}{\mathop{\bf card}}

\newcommand{\strategy}{\mathcal{S}}
\newcommand{\intvars}{\mathcal{I}}
\newcommand{\tightconstraints}{\mathcal{T}}

\newcommand{\loss}{\mathcal{L}}

\newcommand{\Sec}{Section\;}

\newacronym{LO}{LO}{linear optimization}
\newacronym{CO}{CO}{convex optimization}
\newacronym{QO}{QO}{quadratic optimization}
\newacronym{DPP}{DPP}{disciplined parametric program}
\newacronym{SOCO}{SOCO}{second-order cone optimization}
\newacronym{MIQO}{MIQO}{mixed-integer quadratic optimization}
\newacronym{MICO}{MICO}{mixed-integer convex optimization}
\newacronym{MIO}{MIO}{mixed-integer optimization}
\newacronym{LICQ}{LICQ}{linear independence constraint qualification}
\newacronym{BNB}{B\&B}{branch-and-bound}
\newacronym{MILO}{MILO}{mixed-integer linear optimization}
\newacronym{MINLO}{MINLO}{mixed-integer nonlinear optimization}
\newacronym{SPM}{SPM}{successive projection method}
\newacronym{sBB}{sBB}{spacial branch and bound}
\newacronym{NLO}{NLO}{nonlinear optimization}
\newacronym{PWA}{PWA}{piecewise affine}
\newacronym{SVM}{SVM}{support vector machines}
\newacronym{OCT}{OCT}{optimal classification tree}
\newacronym{OCT-H}{OCT-H}{optimal classification trees with-hyperplanes}
\newacronym{CART}{CART}{classification and regression tree}
\newacronym{NN}{NN}{neural network}
\newacronym{ReLU}{ReLU}{rectified linear unit}
\newacronym{MLOPT}{MLOPT}{machine learning optimizer}
\newacronym{MPC}{MPC}{model predictive control}
\newacronym{ADMM}{ADMM}{alternating direction method of multipliers}
\newacronym{UAV}{UAV}{unmanned aerial vehicles}
\newacronym{flop}{flop}{floating point operation}

\begin{document}

\ifpreprint
\title{\bfseries \sffamily Online Mixed-Integer Optimization in Milliseconds}
\author{Dimitris Bertsimas and Bartolomeo Stellato}
\maketitle

\begin{abstract}

We propose a method to solve online~\gls{MIO} problems at very high speed using machine learning.
By exploiting the repetitive nature of online optimization, we are able to greatly speedup the solution time.
Our approach encodes the optimal solution into a small amount of information denoted as {\em strategy} using the Voice of Optimization framework proposed in~\cite{bertsimas2018}.
In this way the core part of the optimization algorithm becomes a multiclass classification problem which can be solved very quickly.
\reviewChanges{In this work,} we extend that framework to real-time and high-speed applications focusing on parametric~\gls{MIQO}.
We propose an extremely fast online optimization algorithm consisting of a feedforward \gls{NN} evaluation and a linear system solution where the matrix has already been factorized.
Therefore, this online approach does not require any solver nor iterative algorithm.
We show the speed of the proposed method both in terms of total computations required and measured execution time.
We estimate the number of~\glspl{flop} required to completely recover the optimal solution as a function of the problem dimensions. Compared to state-of-the-art \gls{MIO} routines, the online running time of our method is very predictable and can be lower than a single matrix factorization time.
We benchmark our method against the state-of-the-art solver Gurobi obtaining from two to three orders of magnitude speedups on examples from fuel cell energy management, sparse portfolio optimization and motion planning with obstacle avoidance.

\end{abstract}

\fi

\glsresetall


\section{Introduction}

\Gls{MIO} has become a powerful tool for modeling and solving real-world decision making problems; see~\cite{juenger2010}.
While most \gls{MIO} problems are $\mathcal{NP}$-hard and thus considered intractable, we are now able to solve instances with complexity and dimensions that were unthinkable just a decade ago.
In~\cite{bixby2010} the authors analyzed the impressive rate at which the computational power of~\gls{MIO} grew in the last 25 years providing over a trillion times speedups.
This remarkable progress is due to both algorithmic and hardware improvements.
Despite these advances, \gls{MIO} is still considered harder to solve than convex optimization and, therefore, it is more rarely applied to online settings.

Online optimization differs from general optimization by requiring on the one hand computing times strictly within the application limits and on the other hand limited computing resources.
Fortunately, while online optimization problems are not the same between each solve, only some parameters vary and the structure remains unchanged.
For this reason, online optimization falls into the broader class of parametric optimization where we can greatly exploit the repetitive structure of the problem instances.
In particular, there is a significant amount of data that we can reuse from the previous solutions.

In a recent work~\cite{bertsimas2018}, the authors constructed a framework to predict and interpret the optimal solution of parametric optimization problems using machine learning.
By encoding the optimal solution into a small amount of information denoted as {\em strategy}, the authors convert the solution algorithm into a multiclass classification problem.
Using interpretable machine learning predictors such as \glspl{OCT}, Bertsimas and Stellato were able to understand and interpret how the problem parameters affect the optimal solutions.
\reviewChanges{Therefore,} they were able to give optimization a {\em voice} that the practitioner can understand.

In this paper we extend the framework from~\cite{bertsimas2018} to online optimization focusing on speed and real-time applications instead of interpretability.
This allows us to obtain an end-to-end approach to solve mixed-integer optimization problems online without the need of any solver nor linear system factorization.
The online solution is extremely fast and can be carried out less than a millisecond reducing the online computation time by more than two orders of magnitude compared to state-of-the-art algorithms.

\subsection{Contributions}
In this work, by exploiting the structure of \gls{MIQO} problems, we derive a very fast online solution algorithm where the whole optimization is reduced to a \gls{NN} prediction and a single linear system solution.
Even though our approach shares the same framework as~\cite{bertsimas2018}, it is substantially different in the focus and the final algorithm. The focus is primarily speed and online optimization applications and not interpretability as in~\cite{bertsimas2018}. This is why\reviewChanges{, for our predictions,} we use non interpretable, but very fast, methods such as~\glspl{NN}. Furthermore, our final algorithm does not involve any convex optimization problem solution as in~\cite{bertsimas2018}. Instead, we just apply simple matrix-vector multiplications. Our specific contributions include:

\begin{enumerate}
  \item We focus on the class of \gls{MIQO} instead of dealing with general \gls{MICO} as in~\cite{bertsimas2018}. This allows us to replace the final step to recover the solution with a simple linear system solution based on the KKT optimality conditions of the reduced problem. Therefore, the whole procedure does not require any solver to run compared to~\cite{bertsimas2018}.

  \item To reduce the number of strategies in larger examples, we reassign the samples to a lower number of selected strategies so that the average suboptimality and infeasibility do not increase above certain tolerances. We define this step as ``strategy pruning'' and formulate it as a large-scale \gls{MILO}. To provide solutions in reasonable times, we develop an approximation algorithm that \reviewChanges{reassigns} the training samples according to the strategies appearing most often.

  \item In several practical applications of \gls{MIQO}, the KKT matrix of the reduced problem does not change with the parameters. In this work we factorize it offline and cache the factorization for all the possible solution strategies appearing in the data. By doing so, our online solution step becomes a sequence of simple forward-backward substitutions that we can carry out very efficiently. Hence, with the offline factorization, our overall online procedure does not even require a single matrix factorization. Compared to~\cite{bertsimas2018}, this further simplifies the online solution step.

  \item After the algorithm simplifications, we derive the precise complexity of the overall algorithm in terms of \glspl{flop} which does not depend on the problem \reviewChanges{parameter} values. This makes the execution time predictable and reliable compared to \gls{BNB} algorithms which often get stuck in the tree search procedure.

  \item We benchmark our method against state-of-the-art \gls{MIQO} solver Gurobi on sparse portfolio trading, fuel battery management and motion planning examples.
	  Thanks to the strategy pruning, we obtain between few hundreds to less than 10,000 strategies for all the examples. This allows to achieve high quality strategy predictions in terms of suboptimality and infeasibility. In particular, the average suboptimality is comparable to the one from Gurobi heuristics and infeasibility is always within acceptable values for the applications considered. \reviewChanges{Timing comparisons on these benchmarks, } show up to three orders of magnitude speedups compared to both Gurobi global optimizer and Gurobi heuristics. The worst-case solution time of our method is also up to three orders of magnitude smaller than the one obtained with \gls{BNB} schemes, enabling real-time implementations in milliseconds.
\end{enumerate}

\subsection{Outline}%
\label{sub:outline}
The structure of the paper is as follows.
In Section~\ref{sec:related_work}, we review recent work on machine learning for optimization outlining the relationships and limitations of other methods compared to approach presented in this work.
In addition, we outline the recent developments in high-speed online optimization and the limited advances that appeared so far for~\gls{MIO}.
In Section~\ref{sec:voice}, we introduce the Voice of Optimization framework from~\cite{bertsimas2018} for general~\gls{MIO} describing the concept of solution \emph{strategy}.
In Section~\ref{sec:machine_learning}, we describe the strategies selection problem as a multiclass classification problem and propose the \gls{NN} architecture used in the prediction phase. We also introduce a strategy pruning scheme to reduce the number of strategies.
Section~\ref{sec:online_optimization} describes the computation savings that we can obtain with problem with a specific structure such as~\gls{MIQO} and the worst-case complexity in terms of number of \glspl{flop}.
In Section~\ref{sec:mlopt}, we describe the overall algorithm with the implementation details.
Benchmarks comparing our method to the state-of-the-art solver Gurobi examples with realistic data appear in Section~\ref{sec:benchmarks}.

\section{Related Work}%
\label{sec:related_work}

\subsection{Machine learning for optimization}
\reviewChanges{Recently,} the operations research community started to focus on systematic ways to analyze and solve combinatorial optimization problems with the help of machine learning.
For an extensive review on the \reviewChanges{topic,} we refer the reader to~\cite{bengio2018}.

Machine learning has so far helped optimization in two directions.
The first one investigates heuristics to improve solution algorithms.
Iterative routines deal with repeated decisions where the answers are based on expert knowledge and manual tuning.
A common example is branching heuristics in~\gls{BNB} algorithms.
\reviewChanges{In general,} these rules are hand tuned and encoded into the solvers.
However, the hand tuning can be hard and is in general suboptimal, especially with complex decisions such as \gls{BNB} algorithm behavior.
To overcome these limitations\reviewChanges{, in~\cite{khalil2016},} the authors learn the branching rules without the need of expert knowledge showing comparable or even better performance than hand-tuned algorithms.
Other promising results using ExtraTrees to learn branching rules appeared in~\cite{alvarez2017}.
We refer the reader to~\cite{Lodi2017} for a review on the intersection of machine learning and branching.

Another example appears in~\cite{bonami2018} where the authors investigate whether it is faster to solve~\glspl{MIQO} directly or as~\gls{SOCO} problems by linearizing the cost function.
This problem becomes a classification problem that offers an advantage based on previous data compared to how the decision is heuristically made inside off-the-shelf solvers.

The second direction poses combinatorial problems as control tasks that we can analyze under the reinforcement learning framework~\cite{sutton2018}.
This is applicable to problems with multistage decisions such as network problems or knapsack-like problems.
\cite{dai2017} learn the heuristic criteria for stage-wise decisions in problems over graphs.
In other words, they build a greedy heuristic framework, where they learn the node selection policy using a specialized neural network able to process graphs of any size~\cite{dai2016}.
For every node, the authors feed a graph representation of the problem to the network and they receive an action-value pair suggesting the next node to select in the optimal path.

Even though these two directions introduce optimization to the benefits of machine learning and show promising results, they do not consider the parametric nature of problems solved in real-world applications.
We often solve the same problem formulation with slightly varying parameters several times generating a large amount of data describing how the parameters affect the solution.

Recently, this idea was used to devise a sampling scheme to collect all the optimal active sets appearing from the problem parameters in continuous convex optimization~\cite{misra2019}.
While this approach is promising, it evaluates online all the combinations of the collected active sets without predicting the optimal ones using machine learning.
Another related work appeared in~\cite{klauco2019} where the authors warm-start online active set solvers using the predicted active sets from a machine learning framework. However, they do not provide probabilistic guarantees for that method and their sampling scheme is tailored to their specific application of~\gls{QO} for \gls{MPC}.

In our work, we propose a new method that exploits the amount of data generated by parametric problems to solve \gls{MIO} online at high speed.
In particular we study how efficient we can make the computations using a combination of machine learning and optimization.
To the authors \reviewChanges{knowledge,} this is the first time machine learning is used to both reduce and make more consistent the solution time of \gls{MIO} algorithms.

\subsection{Online optimization}
\label{sub:online_opt}
Applications of online optimization span a wide variety of fields including scheduling~\cite{catalao2010}, supply chain management~\cite{you2008}, hybrid model predictive control~\cite{bemporad1999}, signal decoding~\cite{damen2000}.

%
%

\paragraph{Embedded optimization.}
Over the last decade there has been a significant attention from the community for tools for generating custom solvers for online parametric programs.
CVXGEN~\cite{mattingley2012} is a code generation software for parametric \gls{QO} that produces a fast and reliable solver. However, its code size grows dramatically with the problem dimensions and it is not applicable to large problems.
More recently, the OSQP solver~\cite{stellato2017a} showed remarkable performance with a light first-order method greatly exploiting the structure of parametric programs to save computation time.
The OSQP authors also proposed an optimized version for embedded applications in~\cite{banjac2017}.
Other solvers that can exploit the structure of parametric programs in online settings include qpOASES~\cite{ferreau2014} for \gls{QO} and ECOS~\cite{domahidi2013} for \gls{SOCO}.

\paragraph{Parametric~\gls{MIQO}.}
\reviewChanges{All previously mentioned approaches focus on continuous convex problems with no integer variables such as \gls{QO}.}
This is because of two main reasons.
On the one hand, mixed-integer optimization algorithms are far more complicated to implement than convex optimization ones since they feature a massive amount of heuristics and preprocessing.
This is \reviewChanges{why,} there is still a huge gap in performance between \reviewChanges{open-source} and commercial solvers for \gls{MIO}.
On the other hand, for many online \reviewChanges{applications,} the solution time required to solve \gls{MIO} problems is still not compatible with the amount of time allowed.
An example is hybrid \gls{MPC} where depending on the system dynamics, we have to solve \gls{MIO} problems online in fractions of a second.
Explicit hybrid \gls{MPC} tackles this issue by precomputing offline the entire mapping between the parameter space to the optimal solution~\cite{bemporad2002}. However, the memory required for storing such solutions grows exponentially with the problem dimensions and this approach easily becomes intractable.

\paragraph{Suboptimal heuristics.}
Other approaches solve these problems only suboptimally using heuristics to deal with insufficient time to compute the globally optimal solution.
Examples include the Feasibility Pump heuristic~\cite{fischetti2005} that iteratively solves \gls{LO} subproblems and rounds their solutions until it finds a feasible solution
for the original \gls{MIO} problem.
Another heuristic works by integrating the \gls{ADMM}~\cite{diamond2018} with rounding steps to obtain integer feasible solutions for the original problem.
The downside of these heuristics is that they do not exploit the large amount of data that we gain by solving the parametric problems over and over again in online settings.

\paragraph{Warm-starting.}
In order to speedup subsequent solves, several works focus on warm-starting \gls{BNB} algorithms~\cite{gamrath2015, ralphs2006}, \reviewChanges{which, in some cases, can significantly reduce the solution time.
However, there can be three possible reasons for which warm-starting can bring no significant benefits.}
First, previous solutions can be infeasible for the current problem and, therefore, not useful to create bounds to quickly prune branches in the~\gls{BNB} tree.
Second, many commercial solvers apply fast heuristics that can quickly obtain good feasible solutions. In case these solutions are as good or better than the provided one, warm-starting does not bring any benefit; see, \eg,~\cite[{\tt Start} variable attribute]{gurobi}.
Third, in~\gls{BNB} algorithms the vast majority of time is usually spent to prove optimality and we are not able to significantly reduce it with a warm-started solution~\cite[Section 4]{gamrath2015}.
Instead of providing only the previous optimal solution, we can pass the previous \gls{BNB} tree and adapt the nodes according to parameter changes~\cite{marcucci2019warm}.
This technique can sometimes greatly reduce the number of \glspl{QO} solved.
However, it still requires a \gls{BNB} algorithm to complete, which might be too slow in fast real-time settings.
\reviewChanges{Compared to warm-starting approaches, our method does not directly exploit the previous solution to accelerate the algorithm. However, it uses history of several previous solution to learn how it changes with data.}

\paragraph{Value function approximations.}
Parametric~\gls{MIQO} have also been studied in terms of how the optimal cost changes with the parameters, \ie, the value function.
The authors of~\cite{hassanzadeh2014} propose an iterative scheme to dynamically generate points to construct approximations of the value function with applications to stochastic integer and bilevel integer optimization problems.
Constructing value functions to solve stochastic \gls{MIO} has been studied also in~\cite{tavaslioglu2019, trapp2013}.
\reviewChanges{Depending on the structure and convexity of the value function approximation, the resulting problem can have different approximation quality and tractability.
In this work, however, instead of computing a value function to reformulate the objective of our optimization problem, we directly encode the optimal solution as the output of a machine learning predictor. In this way, no matter how complex or nonconvex the predictor is, we can obtain very short computation times.}

\paragraph{Our approach.}
Despite all these efforts in solving~\gls{MIO} online, there is still a significant gap between the solution time and the real-time constraints of many applications.
In this \reviewChanges{work,} we propose an alternative approach that exploits data coming from previous problem solutions with high-performance \gls{MIO} solvers to reduce the online computation time making it possible to apply \gls{MIO} to online problems that were not approachable before.
The \reviewChanges{approach} proposed in this paper has already been applied to control problems in robotics in~\cite{cauligi2020learning} by exploiting the application-specific structure of the constraints.




\section{The Voice of Optimization}
\label{sec:voice}

In this section we introduce the idea of the optimal strategy following the same framework first introduced in~\cite{bertsimas2018}.
Given a parameter $\theta \in \reals^p$, we define a {\em strategy} $s(\theta)$ as the complete information needed to efficiently recover the optimal solution of an optimization problem.

Consider the mixed-integer optimization problem
\begin{equation}\label{eq:original_problem}
\begin{array}{ll}
\text{minimize} & f(\theta, x)\\
\text{subject to} & g(\theta, x) \le 0,\\
& x_{\mathcal{I}} \in \integers^d,
\end{array}
\end{equation}
where $x\in\reals^n$ is the decision variable and $\theta \in \reals^p$ defines the parameters affecting the problem.
We denote the cost as $f:\reals^p \times \reals^n \to \reals$ and the constraints as $g: \reals^p \times \reals^n \to \reals^m$.
The vector $x^\star(\theta)$ denotes the optimal solution and $f(\theta, x^\star(\theta))$ the optimal cost function value given the parameter~$\theta$.

\paragraph{Optimal strategy.}
We now define the optimal strategy as the set of tight constraints together with the value of integer variables at the optimal solution.
We denote the {\em tight constraints} $\tightconstraints(\theta)$ as constraints that are equalities at optimality,
\begin{equation}\label{eq:tight_constraints}
	\tightconstraints(\theta) = \{i \in \{1, \dots, m\} \mid g_i(\theta, x^\star(\theta)) = 0\}.
\end{equation}
Hence, given the $\tightconstraints(\theta)$ all the other constraints are redundant for the original problem.

If we assume~\gls{LICQ}, the number of tight constraints is at most $n$ because the tight constraints gradients at the solution $\nabla g_i(\theta, x^\star(\theta)) \in \reals^n,\; i \in \tightconstraints(\theta)$ are linearly independent~\cite[Section 12.2]{nocedal2006}.
When \gls{LICQ} does not hold, the number of tight constraints can be more than the number of decision variables, \ie, $|\tightconstraints(\theta)| > n$.
However, in practice the number of tight constraints is significantly lower than the number of constraints $m$, even for degenerate problems.
This means that for many applications where the number of constraints is larger than the number of variables, by knowing \reviewChanges{$\tightconstraints(\theta)$,} we can neglect a large number of redundant constraints at the optimal solution.

When some of the components of $x$ are integer, we cannot easily compute the optimal solution by knowing only~$\tightconstraints(\theta)$.
This is because the solution retrieval would involve a~\gls{MIO} problem to identify the integer components of $x$.
However, after fixing the integer components to their optimal values $x_{\intvars}^\star(\theta)$, the tight constraints allow us to efficiently compute the optimal solution.
Hence the strategy identifying the optimal solution is a tuple containing the index of tight constraints at optimality and the optimal value of the integer variables, \ie, $s(\theta) = (\tightconstraints(\theta), x_{\intvars}^\star(\theta))$.

\paragraph{Solution method.}
Given the optimal strategy, solving~\eqref{eq:original_problem} corresponds to solving the following optimization problem
\begin{equation}\label{eq:reduced_problem}
\begin{array}{ll}
\text{minimize} & f(\theta, x)\\
\text{subject to} & g_i(\theta, x) \le 0,\quad \forall i \in \tightconstraints(\theta)\\
& x_{\mathcal{I}} = x_{\intvars}^\star(\theta),
\end{array}
\end{equation}
Solving~\eqref{eq:reduced_problem} is much easier than~\eqref{eq:original_problem} because it is continuous, convex and has a smaller number of constraints.
Note that, we cannot in general enforce $g_i(\theta, x) = 0$ for the tight constraints because it would make~\eqref{eq:reduced_problem} nonconvex.  However, we can enforce equalities when $g_i$ are linear in $x$~\cite{boyd2004}.
Since it is a very fast and direct step, we will denote it as {\em solution decoding}.
In Section~\ref{sec:online_optimization} we describe the details of how to exploit the structure of~\eqref{eq:reduced_problem} and compute the optimal solution online at very high speeds.



\section{Machine Learning}
\label{sec:machine_learning}

In this section, we describe how we learn the mapping from the parameters $\theta$ to the optimal strategies $s(\theta)$. In this way, we can replace the hardest part of the optimization routine by a prediction step in a multiclass classification problem where each strategy is a class label.

\subsection{Multiclass Classifier}
\label{sub:multiclass_classifier}

Our classification problem features data points $(\theta_i, s_i),\;i=1,\dots,N$ where $\theta_i \in \reals^p$ are the parameters and $s_i \in \strategy$ the corresponding labels identifying the optimal strategies.
Set $\strategy$ is the set of strategies of cardinality $|\strategy| = M$.
Our goal is to predict $\hat{s}_i$ so that it is as close as possible to the true $s_i$ given sample $\theta_i$.

We solve the classification task using \glspl{NN}.
\glspl{NN} have recently become the most popular machine learning method radically changing the way we classify in fields such as computer vision~\cite{krizhevsky2012}, speech recognition~\cite{hinton2012}, autonomous driving~\cite{bojarski2016}, and reinforcement learning~\cite{silver2017}.

In this work we choose feedforward neural networks because they offer a good balance between simplicity and accuracy without the need of more advanced architectures such as convolutional or recurrent \glspl{NN}~\cite{lecun2015}.

\paragraph{Architecture.}
We deploy a similar architecture as in~\cite{bertsimas2018} consisting of $L$ layers defining a composition of functions of the form
\begin{equation*}
	\hat{s} = h_L(h_{L-1}(\dots h_1(\theta))),
\end{equation*}
where each layer consists of
\begin{equation}\label{eq:nn_layer}
	y_{l} = h_l(y_{l-1}) = \sigma_l(W_l y_{l-1} + b_l), \quad i=1,\dots,L,
\end{equation}
where $y_l \in \reals^{n_l}$.
We define the {\em input} layer as $l=1$ and the {\em output} layer as $l=L$ so that $y_{0} = \theta$ and $y_{L} = \hat{s}$.

Each layer performs an affine transformation with parameters $W_{l}  \in \reals^{n_l \times n_{l-1}}$ and $b_l \in \reals^{n_l}$.
In addition, it includes an activation function $\sigma: \reals^{n_l} \to \reals^{n_l}$ to model nonlinearities.
Inner layers feature a~\gls{ReLU} defined as
\begin{equation*}
	\sigma_l(x) = \max (x, 0),\quad l=1,\dots,L-1,
\end{equation*}
where the $\max$ operator is intended elementwise.
The~\gls{ReLU} operator has become popular because it promotes sparsity to the model outputting $0$ for the negative components of $x$  and because it does not experience vanishing gradient issues typical in sigmoid functions~\cite{goodfellow2016}.

The output layer features a softmax activation function $\sigma_L(x) \in \reals^M$ to provide a normalized ranking between the strategies and quantify how likely they are to be the correct one.
Softmax activation functions are very common in multiclass \reviewChanges{classification} because of their smoothness and the probabilistic interpretation of their output as the relative importance between classes~\cite[Section 6.2.2.3]{goodfellow2016}.
We can write the output layer as
\begin{equation*}
	(\sigma_L(x))_j = \frac{e^{x_j}}{\sum_{j=1}^{M} e^{x_j}},
\end{equation*}
where $0 \le \sigma_L(x) \le 1$ because of the nonnegativity of the exponential functon and the normalization factor.

\paragraph{Learning.}
\reviewChanges{
In order to define a proper cost function and train the network we rewrite the labels as a one-hot encoding, \ie, $s^{\rm oh}_i \in \reals^M$ where $M$ is the total number of classes and all the elements of $s^{\rm oh}_i$ are 0 except the one corresponding to the class which is $1$.
Then we define a smooth cost function, \ie, the cross-entropy loss for which gradient descent-like algorithms work well
\begin{equation*}
	\loss_{\rm NN} =  \sum_{i=1}^{N} -(s^{\rm oh}_i)^\tpose \log(\hat{s}_i),
\end{equation*}
}
where $\log$ is intended elementwise.
This loss $\loss$ can also be interpreted as the distance between the predicted probability density of the labels and to the true one.
The training step consists of applying the classic stochastic gradient descent with the derivatives of the cost function obtained using the back-propagation rule.

\paragraph{Online predictions.}
\reviewChanges{After we complete the model training, we aim at predicting the optimal strategy given $\theta$.}
\reviewChanges{In general, the \gls{NN} prediction works well when the neural network is able to capture the structure of the solution for the values of $\theta$ encounceted in practice.}
However, the model can never be perfect and there can be situations where the prediction is not correct\reviewChanges{, thereby providing suboptimal or infeasible solutions.}
To overcome these possible limitations, instead of considering only the best class predicted by \reviewChanges{the~\gls{NN},} we pick the $k$ most-likely classes.
Afterwards, we can evaluate in parallel their feasibility and objective value and pick the feasible one with the best objective.

\subsection{Strategies Exploration}
\label{sec:strategies_exploration}
It is difficult to estimate the amount of data required to accurately learn the classifier for problem~\eqref{eq:original_problem}.
In particular, given $N$ independent samples $\Theta_N = \{\theta_1,\dots,\theta_N\}$ drawn from an unknown discrete \reviewChanges{distribution,} we find $M$ different strategies $\strategy(\Theta_N)= \{s_1,\dots,s_M\}$. How likely is it to encounter new strategies in the next sample $\theta_{N+1}$?

\reviewChanges{
Following the approach in~\cite{bertsimas2018}, we use the same algorithm to iteratively sample and estimate the probability of finding unseen strategies
\begin{equation*}
	\prob(s(\theta_{N+1}) \notin \strategy(\Theta_N)).
\end{equation*}
First, we define the Good-Turing estimator~\cite{good1953} as
\begin{equation}\label{eq:good_turing}
G = N_1 / N,
\end{equation}
where $N_1$ is the number of distinct strategies appeared exactly once and $N$ the total number of samples.
Then, using the estimator, we bound the probability of finding unseen strategies by applying the following result
\begin{theorem}[Missing strategies bound~\cite{bertsimas2018}]
The probability of encountering a parameter $\theta_{N+1}$ corresponding to an unseen strategy $s(\theta_{N+1})$ satisfies with confidence at least $1 - \beta$
\begin{equation}
\prob(s(\theta_{N+1}) \notin \strategy(\Theta_N)) \le G + c\sqrt{(1/N)\ln(3/\beta)},
\end{equation}
where $G$ corresponds to the Good-Turing estimator~\eqref{eq:good_turing} and $c=(2\sqrt{2} + \sqrt{3})$.
\end{theorem}
In the offline phase, given a desired probability guarantee $\epsilon > 0$ and confidence interval $\beta>0$, we sample strategies and update $G$ until the right hand side bound falls below $\epsilon$, as outline ind Algorithm~\ref{alg:strategies_exploration}~\cite{bertsimas2018}.
\begin{algorithm}
	\caption{Strategies exploration~\cite{bertsimas2018}}
  \label{alg:strategies_exploration}
  \begin{algorithmic}[1]
  \State {\bf given} $\epsilon, \beta, \Theta = \emptyset, \strategy = \emptyset, u = \infty$
  \For{$k = 1,\dots,$}
  \State Sample $\theta_k$ and compute $s(\theta_k)$ \Comment{Sample parameter and strategy.}
  \State $\Theta \gets \Theta \cup \{\theta_k\}$  \Comment{Update set of samples.}
  \If{$s(\theta_k) \notin \strategy$}
  \State $\strategy \gets \strategy \cup \{s(\theta_k)\}$ \Comment{Update strategy set if new strategy found}
  \EndIf
  \If{$G + c \sqrt{(1/k)\ln(3/\beta)} \le \epsilon$} \Comment{Break if bound less than $\epsilon$}
  \State {\bf break}
  \EndIf
  \EndFor
  \State \Return $k, \Theta, \strategy$
  \end{algorithmic}
\end{algorithm}
}

\subsection{Strategy Pruning}%
\label{sub:strategy_pruning}
For some problems, the number of strategies $M$ can quickly grow, thereby making the multiclass classification task very difficult.
Fortunately, different strategies are often redundant because they correspond to multiple global optima and we can select only the relevant strategies to apply.

\paragraph{\Glsdesc{MILO} modelling.}
For every sample $i$ and strategy $j$, we can compute the solution to the reduced problem~\eqref{eq:reduced_problem} obtaining an objective value \reviewChanges{$F_{ij}$}. If the reduced problem is infeasible, we set \reviewChanges{$F_{ij}= \infty$}. To simplify the notation, we refer to \reviewChanges{$f^{\star}_i = f(\theta_i, x^\star(\theta_i))$} as the optimal objective value for sample~$\theta_i$.
We model the sample $i$ to strategy $j$ assignments with variables $x_{ij} \in \{0, 1\}$. Variables $p_j\in\{0, 1\}$ describe \reviewChanges{whether strategy $j$ is picked for discarded.}
The objective is to minimize the number of strategies selected such that the relative cost degradation for each sample $i$ is less than a tolerance $\epsilon$.
The pruning problem can be formulated as the following \gls{MILO},
\reviewChanges{\begin{equation}\label{eq:pruning_problem}
	\begin{array}{ll}
		\text{minimize} & \sum_{j=1}^{M} p_j\\[.5em]
		\text{subject to} & \sum_{j = 1}^{M} F_{ij} z_{ij} \le f^{\star}_i + \epsilon \left|f^{\star}_i\right|,\quad i=1,\dots,N\\[.5em]
		  & \sum_{j=1}^{M} z_{ij} = 1,\quad i=1,\dots,N,\\
				  & z_{ij} \le p_j,\quad i=1,\dots,N,\;j=1,\dots,M,\\
				  & p \in \{0, 1\}^{M},\quad z \in \{0, 1\}^{N \times M},
	\end{array}
\end{equation}}
where $\left|\cdot \right|$ is the absolute value.
Unfortunately, it is very costly to construct problem~\eqref{eq:pruning_problem} because it involves computing $F_{ij}$ for every combination of samples and strategies. For example, if we have 100,000 samples and 2,000 strategies, we need to solve 200,000,000 reduced problems which can be very challenging, even in the specialized cases from Section~\ref{sec:online_optimization}.
In addition, despite recent advances of \gls{MILO} solvers,~\eqref{eq:pruning_problem} with millions of binary variables are often intractable. Therefore, we implement a simpler pruning technique.

\paragraph{Frequency-based heuristic.}
In most cases, the majority of the samples is assigned to a few strategies and the rest of the strategies \reviewChanges{appear} very rarely.
Therefore, if we select the most frequent strategies, we cover the majority of the samples. In this \reviewChanges{way,} we have to reassign only a small portion of the samples to the selected strategies without having to compute $F_{ij}$ for every sample-strategy combination.
The whole procedure is outlined in Algorithm~\ref{alg:strategy_pruning}.
\reviewChanges{Given $\alpha \in (0, 1)$}, the function {\sc SelectFrequentStrategies} in Algorithm~\ref{alg:select_frequent_strategies} selects the most frequent strategies $\strategy_\alpha$ appearing in at least $1- \alpha$ fraction of samples.
\begin{algorithm}
	\caption{({\sc SelectFrequentStrategies}) Select most frequent strategies that are assigned to at least $1-\alpha$ fraction of samples.}
	\label{alg:select_frequent_strategies}
	\begin{algorithmic}
		\State {\bf input }$\alpha, \{s(\theta_i)\}_{i=1}^{N}, \strategy$
		\State {\bf output } $\strategy_{\alpha}$
		\State $t \gets 0,\quad \strategy_{\alpha} \gets \emptyset$
		\For{$s \in \strategy$}
		\State $q_{s} \gets \left|\left\{s(\theta_i) = s,\quad i=1,\dots,N\right\}\right|$ \Comment{Compute strategy occurrences.}
		\EndFor
		\State $v \gets \Call{ReverseArgsort}{q}$ \Comment{Sort strategies by decr.\ occurrences.}
		\For{\reviewChanges{$\ell \in v$}} \Comment{Iterate for every strategy $\ell$}
		\State \reviewChanges{$t \gets t + q_\ell$} \Comment{Update number of samples.}
		\State $\strategy_{\alpha} \gets \strategy_{\alpha} \cup \{\ell\}$
		\IfThen{$t > \lceil (1 - \alpha) N \rceil$}{{\bf break}}
		\EndFor
	\end{algorithmic}
\end{algorithm}
Then, the algorithm selects the discarded samples~$\Theta_d$ that were not assigned to any strategy in $\strategy_{\alpha}$.
These samples are then reassigned by comparing their cost $f_{ij}$ with every selected strategy in $\strategy_{\alpha}$.
If there is at least a sample for which the best reassigned strategy cost \reviewChanges{$r_i$ is above the tolerance, \ie, $r_i > f^{\star}_i + \epsilon \left|f^{\star}_i\right|$}, then $\alpha$ is reduced to $\alpha/2$ to account for more strategies and fewer discarded \reviewChanges{samples} and the iterations are repeated.
Otherwise, we found a strategy pruning and sample assignment satisfying tolerance $\epsilon$ and the algorithm terminates.
If we reach the maximum number of iterations, it means that there is no feasible assignment given the specified tolerance and the last value of fraction of discarded samples $\alpha$.
\begin{algorithm}
	\caption{Prune strategies while keeping feasibility and low suboptimality}
	\label{alg:strategy_pruning}
	\begin{algorithmic}
		\State {\bf input } $\{(\theta_i, s(\theta_i))\}_{1}^{N}, \strategy$
		\State {\bf output } $\strategy$
		\State $\alpha \gets 0.05$
		\For{${\rm it} = 1,\dots, \max_{\rm it}$}
		\State $\strategy_{\alpha} \gets$ \Call{SelectFrequentStrategies}{$\alpha, \{s(\theta_i)\}_{i=1}^{N}, \strategy$}
		\State $\Theta_d \gets \{\theta_i \mid s(\theta_i) \notin \reviewChanges{\strategy_{\alpha}}\}$  \Comment{Select discarded samples}
		\For{$\theta_i \in \Theta_d$}
		\For{$s_j \in \strategy_{\alpha}$}
		\State $F_{ij} \gets \text{Solve}~\eqref{eq:reduced_problem}$ \Comment{Compute sample-strategy pairs}
		\EndFor
		\State \reviewChanges{$r_i \gets \min_j(F_{ij})$} \Comment{Reassign sample $i$ to best strategy}
		\If{\reviewChanges{$r_i \le f^{\star}_i + \epsilon\left|f^{\star}_i\right|$}} \Comment{Suboptimality condition satisfied}
		\State {\bf break}
		\EndIf
		\EndFor
		\State $\alpha \gets \alpha / 2$
		\EndFor
		\State \Return $\strategy$
	\end{algorithmic}
\end{algorithm}
Compared to solving problem~\eqref{eq:pruning_problem}, this method relies on the samples-strategy computation of just a small portion of discarded samples.
In addition, there is no need to solve any large-scale \gls{MILO}, which makes it much more scalable to large settings.
The downside of this heuristic approach is that the best strategy assignment uses the most frequent strategies, which might not always be the optimal one.
However, for small problems where the \gls{MILO} is solvable, the heuristic solution always gives similar number of pruned strategies as \gls{MILO} while always satisfying, by construction, the cost function degradation constraint.

\section{High-Speed Online Optimization}
\label{sec:online_optimization}

Thanks to the learned predictor, our method offers great computational savings compared to solving each problem instance from scratch.
In addition, when the problem offers a specific structure, we can gain even further speedups and replace the whole optimizer with a linear system solution.

\paragraph{Two major challenges in optimization.}
By using our previous solutions, the learned predictor maps new parameters $\theta$ to the optimal strategy replacing two of the arguably hardest tasks in numerical optimization algorithms:
\begin{description}
	\item[Tight constraints.] Identifying the tight constraints at the optimal solution is in general a hard task because of the combinatorial complexity to search over all the possible combinations of constraints. For this \reviewChanges{reason,} the worst-case \reviewChanges{complexity of} active-set methods (also called simplex methods for~\gls{LO}) is exponential in the number of constraints~\cite{bertsimas1997}.

\item[Integer variables.] It is well known that finding the optimal solution of mixed-integer programs is~\NPhard~\cite{bertsimas2005}. Hence, solving problem~\eqref{eq:original_problem} online might require a prohibitive computation time because of the combinatorial complexity of computing the optimal values of the integer variables.
\end{description}
We solve both these issues by evaluating our predictor that outputs the optimal strategy $s(\theta)$, \ie, the tight constraints at optimality $\tightconstraints(\theta)$ and the value of the integer variables $x_{\intvars}^\star(\theta)$.
After evaluating the predictor, computing the optimal solution consists of solving~\eqref{eq:reduced_problem} which we can achieve much more efficiently than solving~\eqref{eq:original_problem}, especially in case of special structure.

\paragraph{Special structure.}

When $g$ is a linear function in $x$ we can directly consider the tight constraints as equalities in~\eqref{eq:reduced_problem} without losing convexity.
In these cases~\eqref{eq:reduced_problem} becomes a convex equality constrained problem that can be solved via Newton's method~\cite[\Sec 10.2]{boyd2004}.
We can further simplify the online solution in special cases such as~\gls{MIQO} (and also \gls{MILO}) of the form
\begin{equation}\label{eq:miqp_online}
\begin{array}{ll}
\text{minimize} & (1/2)x^\tpose P x + q^\tpose x + r\\
\text{subject to} & A x \le b\\
&  x_\intvars \in \integers^d,
\end{array}
\end{equation}
with cost $P \in \symm_{+}^n$, $q\in \reals^n$, $r \in \reals$ and constraints $A \in \reals^{m \times n}$ and $b\in \reals^m$.
We omitted the dependency of the problem data on $\theta$ for ease of notation.
Given the optimal strategy $s(\theta) = (\tightconstraints(\theta), x^\star_{\intvars}(\theta))$, computing the optimal solution to~\eqref{eq:miqp_online} corresponds to solving the following reduced problem from~\eqref{eq:reduced_problem}
\begin{equation}\label{eq:miqp_reduced}
\begin{array}{ll}
\text{minimize} & (1/2)x^\tpose P x + q^\tpose x + r\\
\text{subject to} & A_{\tightconstraints(\theta)} x = b_{\tightconstraints(\theta)}\\
&  x_\intvars = x_\intvars^\star(\theta).
\end{array}
\end{equation}
Since it is an equality constrained \gls{QO}, we can compute the optimal solution by solving the linear system defined by its KKT conditions~\cite[\Sec 10.1.1]{boyd2004}
\begin{equation}
    \begin{aligned}
	\label{eq:kkt_online}
	\begin{bmatrix}
	P & A_{\tightconstraints(\theta)}^\tpose & I_\mathcal{I}^\tpose\\
	A_{\tightconstraints(\theta)} & \multicolumn{2}{c}{\multirow{2}{*}{0}}\\
	\identity_{\intvars} & &
	\end{bmatrix}
	\begin{bmatrix}
	x\\\nu
	\end{bmatrix}
	&= \begin{bmatrix}
	-q\\
	b_{\tightconstraints(\theta)}\\
	x^\star_{\intvars}(\theta)
	\end{bmatrix}.
    \end{aligned}
\end{equation}
Matrix $I$ is the identity matrix.
Vectors or matrices with a subscript index set identify only the rows corresponding to the indices in that set.
$\nu$ are the dual variables of the reduced continuous problem.
The dimensions of the KKT matrix are $q \times q$ where
\begin{equation} \label{eq:kkt_dimension}
q = n + |\tightconstraints(\theta)| + d.
\end{equation}
We can apply the same method to~\gls{MILO} by setting $P=0$.
In case no integer variables are present ($d=0$), the strategy identifies only the tight constraints and the dimension of the linear system~\eqref{eq:kkt_online} reduces to $n + |\tightconstraints(\theta)|$.

\subsection{The Efficient Solution Computation}
\label{sub:efficient_solution_computation}
In case of~\gls{MIQO}, solving the linear system~\eqref{eq:kkt_online} corresponds to computing the solution to~\eqref{eq:miqp_online}.
Let us analyze the components involved in the online computations to further optimize the solution time.

\paragraph{Linear system solution.}
The linear system~\eqref{eq:kkt_online} is sparse and symmetric and we can solve it with both direct methods and indirect methods.
Regarding direct methods, we compute a sparse permuted $LDL^\tpose$ factorization~\cite{davis2006} of the KKT matrix where $L$ is a square lower triangular matrix and $D$ a diagonal matrix both with dimensions $q \times q$.
The factorization step requires $O(q^3)$ number of \glspl{flop} which can be expensive for large systems.
After the factorization, the solution consists only in forward-backward solves which can be evaluated very efficiently and in parallel with complexity $O(q^2)$.
Alternatively, when the system is very \reviewChanges{large,} we can use an indirect method such as MINRES~\cite{paige1982} to iteratively approximate the solution by using simple matrix-vector multiplications at each step with complexity \reviewChanges{$O(q^2)$.}
Note that indirect methods while more amenable for large problems, can suffer from bad scaling of the matrix data, requiring many steps before convergence.

\paragraph{Matrix caching.}
In several \reviewChanges{cases,} $\theta$ does not affect the matrices $P$ and $A$ in~\eqref{eq:miqp_online}.
In other words, $\theta$ enters only in the linear part of the cost and the right hand side of the constraints and does not affect the KKT matrix in~\eqref{eq:kkt_online}.
This means that, since we know all the strategies that appeared in the training phase, we can factor each of the KKT matrices and store the factors $L$ and $D$ offline.
Therefore, whenever we predict the strategy related to a new parameter $\theta$, we can just perform forward-backward solves to obtain the optimal solution without having to perform a new factorization.
This step requires $O(q^2)$~\glspl{flop} -- an order of magnitude less than factorizing the matrix.

\paragraph{Parallel strategy evaluation.}
In Section~\ref{sub:multiclass_classifier} we explained how we transform the strategy selection into a multiclass classification problem.
Since we have the ability to compare the quality of different strategies in terms of optimality and feasibility, online we choose the $k$ most likely strategies from the predictor output and we compare their performance.
Since the comparison is completely independent between the candidate strategies, we can perform it in parallel saving additional computation time.

\subsection{Online Complexity}
We now measure the online complexity of the proposed approach in terms of \glspl{flop}.
The first step consists in evaluating the neural network with $L$ layers.
As shown in~\eqref{eq:nn_layer}, each layer consists in a matrix-vector multiplication and additions $W_l y_{l-1} + b_l$ which has order $O(n_{l}n_{l-1})$ operations~\cite[Section C.1.2]{boyd2004}.
The~\gls{ReLU} step in the layer does not involve any~\gls{flop} since it is a simple truncation of the non positive components of the layer input.
Summing these operations over all the layers, the complexity of evaluating the network becomes
\begin{equation*}
O(n_{1}n_{p} + n_{2}n_{1} + \dots + n_{M}n_{L-1}) = O(n_M n_{L-1}),
\end{equation*}
where $p$ is the dimension of the parameter $\theta$ and we assume that the number of strategies $M$ is larger than the dimension of any layer.
Note that the final softmax layer, while being very useful in the training phase and in its interpretation as likelihood for each class, is unnecessary in the online evaluation since we just need to rank the best strategies.

After the neural network \reviewChanges{evaluation,} we can decode the optimal solution by solving the KKT linear system~\eqref{eq:kkt_online} of dimension defined in~\eqref{eq:kkt_dimension}.
Since we already factorized it, the online computations are just simple forward-backward substitutions as discussed in Section~\ref{sub:efficient_solution_computation}.
Therefore, the \glspl{flop} required to solve the linear system are
\begin{equation*}
O((n + |\tightconstraints(\theta)| + d)^2).
\end{equation*}
This dimension is in general much smaller than the number of constraints $m$ and mostly depends only on the problem variables.
In addition the KKT matrix~\eqref{eq:kkt_online} is sparse and if the factorization matrices are sparse as well, we could further reduce the \glspl{flop} required.

The overall complexity of the complete \gls{MIQO} online solution taking into account both \gls{NN} prediction and solution decoding becomes
\begin{equation}
  O(n_M n_{L-1} + (n + |\tightconstraints(\theta)| + d)^2),
\end{equation}
which does not depend on the cube of any of the problem dimensions.
This means that with dense matrices, our method is asymptotically cheaper than factorizing a single linear system since that operation would scale with the cube of its input data.
For sparse matrices, we can make similar considerations based on the number of nonzero elements instead of the dimensions~\cite{davis2006}.

Despite these high-speed considerations on the number of operations required, it is very important to remark the reliability of the computation time required by our approach.
The execution time of \gls{BNB} algorithms greatly depends on how the solution search tree is analyzed and pruned.
This can vary significantly when problem parameters change making the whole online optimization procedure unreliable in real-time applications.
On the contrary, our method offers a fixed number of operations which we evaluate every time we encounter a new parameter.



\section{Machine Learning Optimizer}
\label{sec:mlopt}

Our implementation extends the software tool \gls{MLOPT} from~\cite{bertsimas2018} which is implemented in Python and integrated with CVXPY~\cite{diamond2016} to model the optimization problems\gls{MLOPT} is available at
\begin{center}
	\url{https://github.com/bstellato/mlopt}.
\end{center}
To speedup the repetitive canonicalizations of parametric~\gls{MIQO} we use the~\gls{DPP} language introduced in CVXPY 1.1~\cite{cvxpylayers2019}. In this way, independently from how the parameters affect the data, constructing a new problem instance given a new $\theta_i$ \reviewChanges{consists of a sparse matrix-vector multiplication~\cite[Section 4.2, matrix $C$]{cvxpylayers2019}}.
MLOPT relies on the Gurobi Optimizer~\cite{gurobi} to solve the problems in the training phase with high accuracy and identify the tight constraints.

After \gls{MLOPT} collects the strategies, we apply Algorithm~\ref{alg:strategy_pruning} to select the most frequent ones and reassign the samples accordingly. This process is performed in parallel over multiple cores to minimize the independent evaluations.
Then, MLOPT passes the data to PyTorch~\cite{paszke2017} using Pytorch-Lightning~\cite{falcon2019pytorch} library to define the architecture and train the \gls{NN} to classify the strategies.
We split the training data into $80\;\%$ training and $20\;\%$ validation.
\reviewChanges{
	We tune the neural network parameters (interval): depth ($[3, 15]$), width ($[4, 128]$), learning rate ($[10^{-5}, 10^{-1}]$), batch size ($[32, 256]$), and number of epochs ($[5, 30]$) using the Optuna hyperparameter framework~\cite{optuna_2019} to exploit high parallelization and early pruning.
}

In addition to the \gls{MLOPT} framework in~\cite{bertsimas2018}, we include a specialized solution method for \gls{MIQO} based on the techniques described in Section~\ref{sec:online_optimization} where we factorize and cache the factorization of the KKT matrices~\eqref{eq:kkt_online} for each strategy of the problem to speedup the online computations. \reviewChanges{Note that the memory requirements of this step are not limiting since the stored matrices are in general very sparse and involve a reduced version of the original optimization problem.}

As outlined in Section~\ref{sec:online_optimization}, when we classify online using the \reviewChanges{\gls{NN},} we pick the best $k$ strategies and evaluate them in parallel, in terms of objective function value and infeasibility, to choose the best one.
This step requires a minimal overhead since the matrices are all factorized and we can execute the evaluations in parallel.

\begin{figure}
  \centering
  {\bfseries \sffamily Offline}\\[1em]
\begin{tikzpicture}[auto,>=latex', node distance=8em]
    \node [block] (cvxpy) {CVXPY\\Modeling};
    \node [block] (sampling) [right of=cvxpy] {Strategies\\Sampling};
    \path[->, draw] (cvxpy) to (sampling);
    \node [block] (training) [right of=sampling] {NN\\Training};
    \path[->, draw] (sampling) to (training);
    \node [block] (caching) [right of=training] {Factorization\\and Caching};
    \path[->, draw] (training) to (caching);
\end{tikzpicture}\\[.5em]
  {\bfseries \sffamily Online}\\[1em]
\begin{tikzpicture}[auto,>=latex', node distance=7em]
    \node [block] (prediction) {Strategy\\Prediction};
    \node (theta)[left of=prediction] {$\theta$};
    \path[->, draw] (theta) to (prediction);
    \node (strategy) [right of=cvxpy] {$s(\theta)$};
    \path[->, draw] (prediction) to (strategy);
    \node [block] (decoding) [right of=strategy] {Solution\\Decoding};
    \path[->, draw] (strategy) to (decoding);
    \node (x)[right of=decoding] {$x^\star$};
    \path[->, draw] (decoding) to (x);
\end{tikzpicture}\\[.5em]
\caption{Algorithm implementation}
\label{fig:block_diagram}
\end{figure}
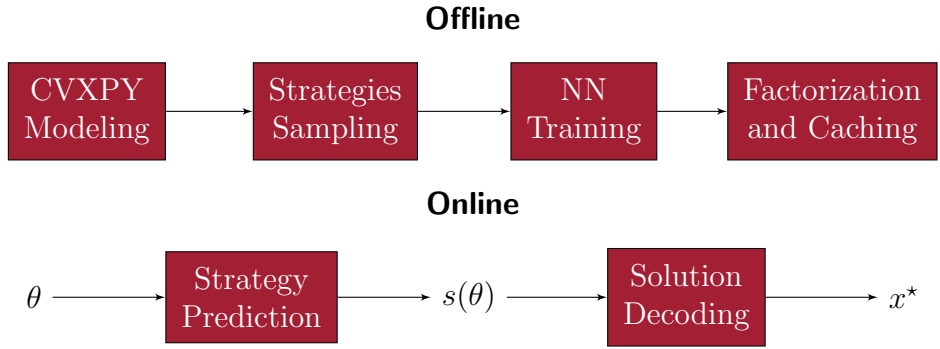

We also parallelize the training phase where we collect data and solve the problems over multiple CPUs. The \gls{NN} training takes place on a GPU which greatly reduces the training time.
An overview of the online and offline algorithms appears in Figure~\ref{fig:block_diagram}.

\section{Computational Benchmarks}
\label{sec:benchmarks}

In this section, we benchmark our learning scheme on multiple parametric examples from continuous and mixed-integer problems.
We compare the predictive performance and the computational time required by our method to solve the problem compared to using GUROBI Optimizer~\cite{gurobi}.
\reviewChanges{We run Gurobi with warm-starting enabled in order to reuse the solution obtained from the previous parameter value $\theta_i$. We execute it in two variants: default settings and ``heuristic'' mode with time limit of 1 second ({\tt TimeLimit=1}) and focus on feasibility ({\tt MIPFocus=1}).}
\reviewChanges{We report the execution time of all the compared methods in seconds. In these example, the MLOPT \gls{NN} prediction takes always less than $0.5\;ms$. Therefore, almost all the reported MLOPT time consists of strategies evaluation.}
\reviewChanges{We run the experiments on the Princeton Institute for Computational Science and Engineering (PICSciE) facility exploiting 16 parallel Intel Broadwell e5-2680v4 cores for the data collection involving the problems solution and a NVIDIA P100 GPUs for the neural network training.
We execute both MLOPT and GUROBI only on CPUs in the online phase.}
For each \reviewChanges{problem, we sample 100,000} parameters $\theta_i$ to collect the strategies.  We choose this number because the~\gls{NN} training works better with a large number of data points.
Note that, thanks to the multiple cores and the code parallelizations, we were able to train the algorithm between a few hours and less than a day, even for problems that take up to hundreds of seconds to solve with Gurobi.
\reviewChanges{
We run MLOPT with the default parameters described in Section~\ref{sec:mlopt}.
}

In addition, for all the examples, the Good-Turing estimator condition from Section~\ref{sec:strategies_exploration} was always satisfied for $\epsilon = 0.001$.
We use $10,000$ samples in the test set of the algorithm.

\paragraph{Infeasibility and suboptimality.}
We follow the same criteria for calculating the relative suboptimality and infeasibility as in~\cite{bertsimas2018}. We repeat them here for completeness.
After the learning phase, we compare the predicted solution $\hat{x}^\star_i$ to the optimal one $x^\star_i$ obtained by solving the instances from scratch.
Given a parameter $\theta_i$, we say that the predicted solution is infeasible if the constraints are violated more than $\epsilon_{\rm inf} = 10^{-4}$ according to the infeasibility metric
\begin{equation*}
	p(\hat{x}^\star_i) = \|(g(\theta_i, \hat{x}^\star_i))_{+}\|_{\infty} / r(\theta_i, \hat{x}^{\star}_i),
\end{equation*}
where $r(\theta, x)$ normalizes the violation depending on the size of the summands of $g$.
In case of \gls{MIQO}, $g(\theta, x) = A(\theta)x - b(\theta)$ and ${r(\theta, x) = \|b(\theta)\|_{\infty}}$.
If the predicted solution $\hat{x}^\star_i$ is feasible, we define its suboptimality as
\begin{equation*}
	d(\hat{x}^\star_i) = (f(\theta_i, \hat{x}^\star_i) - f(\theta_i, x^\star_i))/|f(\theta_i, x^\star_i)|,
\end{equation*}
where $f$ is the objective of our~\gls{MIQO} problem in~\eqref{eq:miqp_online}. Note that $d(\hat{x}^\star_i) \ge 0$ by construction.
For each \reviewChanges{example}, we report the average infeasibility and suboptimality over the samples.

\reviewChanges{
\paragraph{Accuracy.}
In multiclass classification, accuracy corresponds to the fraction of times the predicted class is correct.
However, in this setting, we can have multiple strategies (classes) leading to high quality solutions.
Therefore, we adapt the concept of accuracy to take into account the quality of the solutions, instead of the specific strategy used during training.
In other words, we consider a predicted solution to be accurate if it is feasible and if the suboptimality is less than the tolerance $\epsilon_{\rm sub} = 10^{-4}$.
We, then, define the accuracy as
\begin{equation*}
(1/N) |\{\hat{x}_i^\star \mid p(\hat{x}_i^\star) \le \epsilon_{\rm inf}\quad \mbox{and} \quad d(\hat{x}_i^\star) \le \epsilon_{\rm sub}\}|.
\end{equation*}

}

\subsection{Fuel Cell Energy Management}%
\label{sub:fuel_cell_problem}

Fuel cells are a green highly efficient energy source that need to be controlled to stay within admissible operating ranges.
Too large switching between ON and OFF states can reduce both the lifespan of the energy cell and increase energy losses.
This is why fuel cells are often paired with energy storage devices such as capacitors which help reducing the switching frequency during fast transients.

In this \reviewChanges{example,} we would like to control the energy balance between a super capacitor and a fuel cell in order to match the demanded power~\cite{frick2015}.
The goal is to minimize the energy losses while maintaining the device within acceptable operating limits to prevent lifespan degradation.

We can model the capacitor dynamics as
\begin{equation}
	E_{t+1} = E_{t} + \tau (P_t - P^{\rm load}_t),
\end{equation}
where $\tau > 0$ is the sampling time and $E_t \in [E^{\rm min}, E^{\rm max}]$ is the energy stored.
$P_t \in [0, P^{\rm max}]$ is the power provided by the fuel cell and $P^{\rm load}$ is the desired load power.

At each time $t$ we model the on-off state of the fuel cell with the binary variable $z_t \in \{0, 1\}$.
When the battery is off ($z_t = 0$) we do not consume any energy, thus we have $0 \le P_{t} \le P^{\rm max} z_t$.
When the engine is on ($z_t = 1$) it consumes $\alpha P_{t}^2 + \beta P_{t} + \gamma$ units of fuel, with $\alpha, \beta, \gamma > 0$.
We define the stage power cost as
\begin{equation*}
    f(P, z) = \alpha P^2 + \beta P + \gamma z.
\end{equation*}

We now model the total sum of the switchings over a time window in order to constrain its value.
In order to do so we introduce binary variable $d_t \in \{0, 1\}$ determining whether the cell switches at time $t$ either from ON to OFF or viceversa.
Additionally we introduce the auxiliary variable $w_t \in [-1, 1]$ accounting for the amount of change brought by $d_t$ in the battery state $z_t$,
\begin{equation}
w_t = \begin{dcases}
1, & d_t = 1 \wedge z_t = 1,\\
-1, & d_t = 1 \wedge z_t = 0,\\
0, & \text{otherwise}.
\end{dcases}
\end{equation}
We can model these logical relationships as the following linear inequality~\cite{frick2015},
\begin{equation}
G (w_t, z_t, d_t) \le h, \quad \text{with}\quad
	G = \begin{bmatrix}
	1 & 0 & -1\\
	-1 & 0 & -1\\
	1 & 2 & 2\\
	-1 & -2 & 2
	\end{bmatrix}, \quad
	h = (0, 0, 3, 1).
\end{equation}
Hence we can write the number of switchings $s_{t+1}$ appeared up to time $t+1$ over the past time window of length $T$ as
\begin{equation}
	s_{t+1} = s_{t} + d_{t} - d_{t - T},
\end{equation}
and impose the constraints $ s_t \le n^{\rm sw}$.
The complete fuel cell problem becomes
\begin{equation}\label{eq:fuel_cell_problem}
\begin{array}{ll}
    \mbox{minimize}   &  \displaystyle \sum_{t=0}^{T - 1} f(P_{t}, z_t)\\[1em]
    \mbox{subject to} & E_{t + 1} = E_{t} + \tau (P_t - P_{t}^{\rm load}),\\
    & E^{\rm min} \le E_t \le E^{\rm max},\\
    & 0 \le P_t \le z_t P^{\rm max},\\
    & z_{t+1} = z_t + w_t,\\
    & s_{t+1} = s_{t} + d_{t} - d_{t - T},\\
    & s_t \le n^{\rm sw},\\
    & G (w_t, z_t, d_t) \le h,\\
    & E_0 = E_{\rm init}, \quad z_0 = z_{\rm init},\quad s_0 = s_{\rm init}\\
    & z_t \in \{0, 1\}, \quad d_t \in \{0, 1\},\quad w \in [-1, 1].
\end{array}
\end{equation}
The problem parameters are $\theta = (E_{\rm init}, z_{\rm init}, s_{\rm init}, d^{\rm past}, P^{\rm load})$ where $d^{\rm past}=(d_{-T}, \dots, d_{-1})$ and $P^{\rm load} = (P^{\rm load}_0, \dots, P^{\rm load}_{T-1})$.
In order to properly control the dynamical system, we must solve~\eqref{eq:fuel_cell_problem} within each sampling time $\tau$.

\paragraph{Problem setup.}
We chose parameters from~\cite{frick2015} with values ${\alpha= 6.7 \times 10^{-4}}$, $\beta = 0.2$, $\gamma = 80\;\si{\watt}$, and sampling time $\tau=1\;\si{second}$
We define energy and power constraints with $E^{\rm min} = 5.2\;\si{\kilo \joule}$, $E^{\rm max} = 10.2\;\si{\kilo \joule}$ and $P^{\rm max} = 1.2\;\si{\kilo \watt}$.
The initial values are $E_{\rm init} = 7.7\;\si{\kilo\joule}$, $z_{\rm init} = 0$ and $s_{\rm init}=0$.
We randomly generated the load profile $P^{\rm load}$.

We obtain the offline samples by simulating a closed-loop trajectory of $10,000$ time steps and storing the parameters for each component of $\theta$ along the trajectory \ie, $E_{\rm init}, z_{\rm init}, s_{\rm init}, d^{\rm past}$, and $P^{\rm load}$.
We, then, sample from a uniform distribution over a hypershpere of radius $0.5$ centered at each trajectory point.
Afterwards, we enforce the feasibility of the problem parameters according to the constraints in~\eqref{eq:fuel_cell_problem}.

\paragraph{Results.}
Table~\ref{tab:fuel_cell} reports the problem dimensions and the maximum computation time needed with each technique.
\reviewChanges{The strategy pruning is able to, sometimes, significantly $M$. For example, for $T=40$ we are able to reduce the number of unique strategies $M$ from $8604$ to $1362$.}
Figure~\ref{fig:fuel_cell_mlopt} shows the performance of MLOPT for varying values of the top-$k$ strategies and in terms of computation time, suboptimality, infeasibility and accuracy.
As expected, as $k$ increases the performance improves.
For time horizons of $T=50$ or $T=60$, Gurobi takes on average longer than the allowed sampling time, $\tau=1$ sec.
As noted by the authors of~\cite{frick2015}, in order to get good performance with this system, horizons in at least of $T=60$ need to be considered and the respective solutions are not computable in less than $\tau=1$ sec with state-of-the-art algorithms.
From Table~\ref{tab:fuel_cell}, the maximum time of Gurobi is above the allowed sampling time for horizons $\ge 20$ and Gurobi heuristic for horizon $\ge 40$. Therefore, these solution methods are not applicable for real-time optimization of this dynamical \reviewChanges{system.}

\begin{figure}
	\begin{center}
	\small
\ref{perfcontrol}

\begin{tikzpicture}
	\begin{axis}
		[
		axis x line=bottom,
		axis y line=left,
		enlarge x limits=upper,  
		enlarge y limits=auto,
		ymode=log,          
		width=.45\textwidth,
		height=0.2\textheight,
		xlabel=$k$,
		ymin=1e-05,
		ymax=1e1,
		ylabel={Computation time [s]},
		xtick = data,
		cycle list name=black white,
		thick,
		legend columns=-1,
		legend to name=perfcontrol,
		legend entries = {$T=10$, $T=20$, $T=30$, $T=40$, $T=50$, $T=60$},
		]
		\foreach \t in {10, 20, 30, 40, 50, 60}{
			\addplot+
				table [x=n_best,
					y=mean_time_pred,
					col sep=comma] {./data/benchmarks/control/control_\t_performance.csv};
		}
     \addplot[mark=none, dashed, thick] coordinates {(1,1.0) (100,1.0)};  
	\end{axis}
\end{tikzpicture}
\begin{tikzpicture}
	\begin{axis}
		[
		axis x line=bottom,
		axis y line=left,
		enlarge x limits=upper,  
		enlarge y limits=auto,
		ymode=log,          
		width=.45\textwidth,
		height=0.2\textheight,
		xlabel=$k$,
		ylabel={Suboptimality},
		ymax=1,
		ymin=1e-05,
		xtick = data,
		cycle list name=black white,
		thick,
		]
		\foreach \t in {10, 20, 30, 40, 50, 60}{
			\addplot+ [error bars/.cd, y explicit, y dir=both]
				table [x=n_best,
					y=avg_subopt,
					col sep=comma] {./data/benchmarks/control/control_\t_performance.csv};
		}
	\end{axis}
\end{tikzpicture}

\begin{tikzpicture}
	\begin{axis}
		[
		axis x line=bottom,
		axis y line=left,
		enlarge x limits=upper,  
		enlarge y limits=auto,
		ymode=log,          
		width=.45\textwidth,
		height=0.2\textheight,
		xlabel=$k$,
		ylabel={Infeasibility},
		ymax=1,
		ymin=1e-06,
		xtick = data,
		cycle list name=black white,
		thick,
		]
		\foreach \t in {10, 20, 30, 40, 50, 60}{
			\addplot+
				table [x=n_best,
					y=avg_infeas,
					col sep=comma] {./data/benchmarks/control/control_\t_performance.csv};
		}
	\end{axis}
\end{tikzpicture}
\begin{tikzpicture}
	\begin{axis}
		[
		axis x line=bottom,
		axis y line=left,
		enlarge x limits=upper,  
		enlarge y limits=auto,
		width=.45\textwidth,
		height=0.2\textheight,
		ymax=110,
		ymin=0,
		xlabel=$k$,
		ylabel={Accuracy [\%]},
		xtick = data,
		cycle list name=black white,
		thick,
		]
		\foreach \t in {10, 20, 30, 40, 50, 60}{
			\addplot table [x=n_best, y=accuracy, col sep=comma] {./data/benchmarks/control/control_\t_performance.csv};
	        }
	\end{axis}
\end{tikzpicture}
\end{center}
\caption{MLOPT average performance indicators for fuel cell battery management example. The dashed line indicates the sampling time.}
\label{fig:fuel_cell_mlopt}
\end{figure}
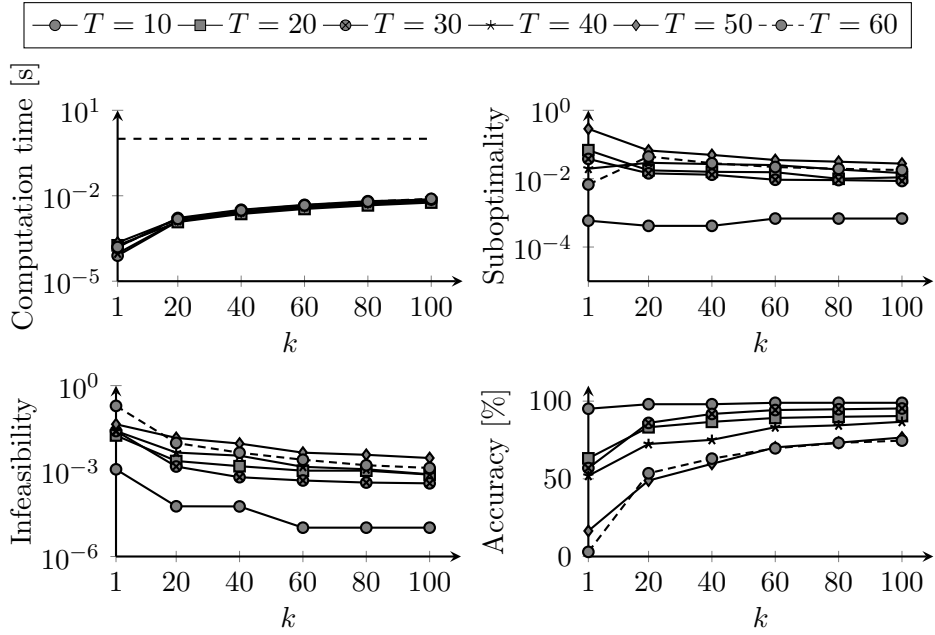

\begin{figure}
	\begin{center}
	\small
\ref{gurobicontrol}

\begin{tikzpicture}
	\begin{axis}
		[
		axis x line=bottom,
		axis y line=left,
		enlarge x limits=upper,  
		enlarge y limits=auto,
		ymode=log,          
		width=.45\textwidth,
		height=0.2\textheight,
		xlabel=$T$,
		ylabel={Computation time [s]},
		ymin=1e-05,
		ymax=1e1,
		xtick = data,
		cycle list name=black white,
		thick,
		legend columns=-1,
		legend to name=gurobicontrol,
		legend entries = {Gurobi heuristic, MLOPT $(k=100)$, Gurobi},
		]
		\addplot+
			table [x=T,
				y=mean_time_heuristic,
				col sep=comma] {./data/benchmarks/control/control_performance_full_heuristic.csv};
		\addplot+ 
			table [x=T,
				y=mean_time_pred,
				col sep=comma] {./data/benchmarks/control/control_performance_full_heuristic.csv};
		\addplot+ 
			table [x=T,
				y=mean_time_full,
				col sep=comma] {./data/benchmarks/control/control_performance_full_heuristic.csv};
     \addplot[mark=none, dashed, thick] coordinates {(10,1.0) (60,1.0)};  
	\end{axis}
\end{tikzpicture}
\begin{tikzpicture}
	\begin{axis}
		[
		axis x line=bottom,
		axis y line=left,
		enlarge x limits=upper,  
		enlarge y limits=auto,
		ymode=log,          
		width=.45\textwidth,
		height=0.2\textheight,
		xtick = {10, 20, 30, 40, 50, 60},
		xlabel=$T$,
		ylabel={Suboptimality},
		ymin=1e-05,
		ymax=1,
		cycle list name=black white,
		thick,
		]
		\addplot+
			table [x=T,
				y=avg_subopt_heuristic,
				col sep=comma] {./data/benchmarks/control/control_performance_full_heuristic.csv};
		\addplot+
			table [x=T,
				y=avg_subopt,
				col sep=comma] {./data/benchmarks/control/control_performance_full_heuristic.csv};
	\end{axis}
\end{tikzpicture}
\end{center}
\caption{Comparison between Gurobi and MLOPT performance for the fuel cell battery management example. Average computation time and suboptimality. The dashed line indicates the sampling time. \reviewChanges{Missing points on the suboptimality plot correspond to machine precision $10^{-15}$.}}
\label{fig:l_gurobi}
\end{figure}
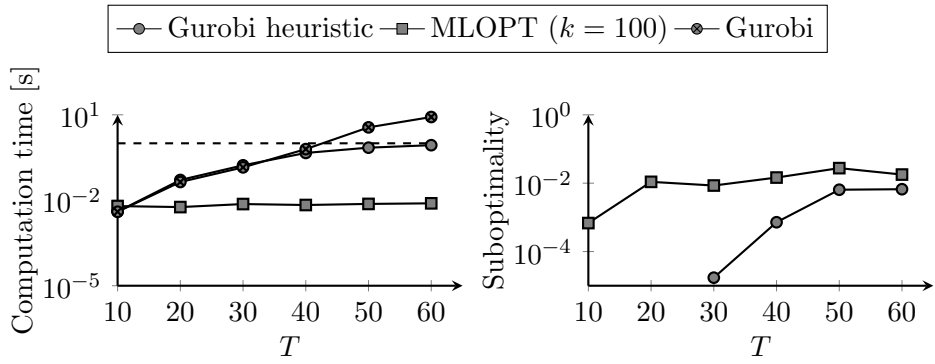

\begin{table}
	\footnotesize
  \centering
  \caption{Fuel cell energy management problem dimensions and maximum times.}
  \label{tab:fuel_cell}
  \begin{tabular}
	  {@{}lllll
    S[round-precision = 4, table-number-alignment=left,scientific-notation=false]
    S[table-format=2.4, round-precision = 4, table-number-alignment=left,scientific-notation=false]
    S[round-precision = 4, table-number-alignment=left,scientific-notation=false]
    @{}}
    \toprule
    $T$ & $n_{\rm var}$ & $n_{\rm constr}$ & {\makecell[lb]{$M$\\ (unpruned)}} & $M$ & {\makecell[lb]{$t_{\rm max}$\\ MLOPT [s]}} & {\makecell[lb]{$t_{\rm max}$\\ Gurobi [s]}}& {\makecell[lb]{$t_{\rm max}$\\ Gurobi \\heuristic [s]}}\\
    \midrule
    \begin{filecontents*}{./data/benchmarks/control/control_performance_full_heuristic.csv}
,problem,learner,n_var,n_constr,n_theta,n_strategies,n_strategies_unpruned,n_infeas,avg_infeas,std_infeas,max_infeas,avg_subopt,std_subopt,max_subopt,avg_subopt_heuristic,std_subopt_heuristic,max_subopt_heuristic,mean_solve_time_pred,std_solve_time_pred,mean_pred_time_pred,std_pred_time_pred,mean_time_pred,std_time_pred,max_time_pred,mean_time_full,std_time_full,max_time_full,mean_time_heuristic,std_time_heuristic,max_time_heuristic,T
5,problem,pytorch,88,207,23,467,622,32,1.0295112739957227e-05,0.0002277258737983,0.013392857142857,0.0006813715422236,0.0086582376663033,0.1499588419959434,-6.383735988947536e-17,1.3002700234041885e-15,0.0,0.0062449586354423,0.0001201058398722,7.461013350598067e-06,1.694065894508601e-21,0.0062524196487929,0.0001201058398722,0.009605239333663,0.0039572610655405,0.0013909521646145,0.0088810920715332,0.0040129309318763,0.0014289309601426,0.0096759796142578,10
5,problem,pytorch,178,417,43,1298,1375,377,0.0007501436335073,0.0052988484998632,0.0645161290322581,0.0109370205385922,0.0661544115817649,2.9523969047616623,-1.2290898212824881e-15,1.3165257696620511e-14,1.28701695727801e-13,0.0057503390822014,0.0002094143480149,2.0558953033897627e-05,3.3881317890172014e-21,0.0057708980352353,0.0002094143480149,0.0095577789542546,0.0437737758494012,0.0729399810388777,0.7645981311798096,0.0515402336817598,0.0742196819914822,0.793212890625,20
5,problem,pytorch,268,627,63,4522,5143,161,0.0003704778519894,0.0048452440130751,0.0892857142857142,0.0085409956901389,0.075378143547844,2.952396904761657,1.7213791884907978e-05,0.0014921895801556,0.1467635877092905,0.0073053320700608,0.0003494524577403,6.421494203342443e-05,1.3552527156068804e-20,0.0073695470120942,0.0003494524577403,0.0107122271032883,0.142932624623109,0.2153704354132885,1.5771291255950928,0.1687566828001026,0.2520720931039161,1.3413901329040527,30
5,problem,pytorch,358,837,83,1362,8604,432,0.0007682295416672,0.0046075735701932,0.0892857142857267,0.0146093552277885,0.1007642225465748,1.676559264706433,0.0007196213166175,0.0199900833581617,0.9426369934449306,0.006801547476414,0.0003874074626769,2.0164766206867253e-05,6.776263578034403e-21,0.0068217122426209,0.0003874074626769,0.0124660914324666,0.6255168079933372,0.7990245885162613,5.469780921936035,0.4588031597693064,0.445158055758245,1.654430866241455,40
5,problem,pytorch,448,1047,103,3955,13590,1088,0.0028688103499436,0.012337784411737,0.1067708333334157,0.0275565657197399,0.1868046596933824,5.166694583328075,0.0064370025662726,0.077915216785735,3.250594044691551,0.0073803274718911,0.0005690225942405,6.624629143307563e-05,0.0,0.0074465737633242,0.0005690225942405,0.013010229369436,3.625787942677529,4.997616847542691,44.93497705459595,0.702877153496395,0.4372885027262622,1.87005615234375,50
5,problem,pytorch,538,1257,123,4734,22691,848,0.001284527395506,0.0057133800665482,0.0747767857142855,0.0178710873222011,0.0693349008112255,2.952396904760173,0.0066596698348786,0.0780941804763255,2.886178974483481,0.0077396957826763,0.0001407080794156,7.507755452111572e-05,0.0,0.0078147733371975,0.0001407080794156,0.0098001713960738,8.366851403927733,12.865871573513054,114.81485199928284,0.8585460191071728,0.3662635675433924,1.6339130401611328,60
\end{filecontents*}
\csvreader[head to column names, late after line=\\]{./data/benchmarks/control/control_performance_full_heuristic.csv}{
    T=\T,
    n_constr=\nconstr,
    n_var=\nvar,
    n_strategies=\M,
    n_strategies_unpruned=\Munpruned,
    max_time_pred=\tmaxmlopt,
    max_time_full=\tmaxgurobi,
    max_time_heuristic=\tmaxheuristic,
    }{\T & \nvar & \nconstr & \Munpruned & \M & \tmaxmlopt & \tmaxgurobi & \tmaxheuristic}
\bottomrule
  \end{tabular}
\end{table}

\subsection{Portfolio Trading}%
\label{sub:portfolio_trading}
Consider the portfolio investment problem~\cite{markowitz1952, boyd2017}.
This problem has been extensively analyzed in robust optimization~\cite{bertsimas2008} and stochastic control settings~\cite{herzog2007}.
The decision variables are the normalized portfolio weights $w_t\in \reals^{n+1}$ at time $t$ corresponding to each of the assets in the portfolio and a riskless asset at the $(n+1)$-th position denoting a cash account.
We define the trade as the difference of consecutive weights $w_{t} - w_{t-1}$ where $w_{t-1}$ is the given vector of current asset investments acting as a problem parameter.
The goal is to maximize the risk-adjusted returns as follows
\begin{equation}\label{eq:multiperiod_portfolio}
\begin{array}{ll}
\text{maximize} & \hat{r}_t^\tpose w_t - \gamma \ell^{\rm risk}_t(w_t)  - \ell^{\rm hold}_t(w_t) - \ell^{\rm trade}_t(w_t - w_{t-1})\\[.5em]
\text{subject to} & \ones^\tpose w_t = 1,\\
& \card(w_t) \le c.
\end{array}
\end{equation}
Four terms compose the stage rewards.
First, we describe returns $\hat{r}_t^\tpose w_t$ as a function of the estimated stock returns $\hat{r}_t \in [0, 1]^{n+1}$ at time $t$.
Second, we define the risk cost as
\begin{equation*}
	\ell^{\rm risk}_t(x) = x^\tpose \hat{\Sigma}_t x,
\end{equation*}
where $\hat{\Sigma}_t \in \symm^{(n+1) \times (n+1)}_{+}$ is the estimate of the covariance of the returns at time $t$.
Third, we define the holding cost as
\begin{equation*}
	\ell^{\rm hold}_t(x) = s_t^\tpose (x)_{-},
\end{equation*}
where $(s_t)_i \ge 0$ is the borrowing fee for shorting asset $i$ at time $t$.
The fourth term describes a penalty on the trades defined as
\begin{equation*}
	\ell^{\rm trade}_t(x) = \lambda\|x_t - x_{t-1}\|_1,
\end{equation*}
Parameter $\lambda > 0$ denotes the relative cost importance of penalizing the trades.
The first constraint enforces the portfolio weights normalization while the second constraints the maximum number of nonzero asset investments, \ie, the cardinality, to be less than $c \in \integers_{> 0}$.
For the complete model derivation without cardinality constraints see~\cite[Section 5.2]{boyd2017}.

\paragraph{Risk model.}
We use a common risk model described as a $k$-factor model $\Sigma_t = F_t \Sigma_t^F F_t^\tpose + D_t$ where $F_t \in \reals^{(n+1) \times k}$ is the factor loading matrix and $\Sigma_t^F \in \symm_{+}^{k \times k}$ is an estimate of the factor returns $F^\tpose r_t$ covariance~\cite{boyd2017}.
Each entry $(F_t)_{ij}$ is the loading of asset $i$ to factor $j$.
$D_t \in \symm_{+}^{(n+1) \times (n+1)}$ is a nonnegative diagonal matrix accounting for the additional variance in the asset returns usually called idiosyncratic risk.
We compute the factor model estimates $\hat{\Sigma}_t$ with $15$ factors by using a similar method as in~\cite[\Sec 7.3]{boyd2017} where we take into account data from the two years time window before time~$t$.

\paragraph{Return forecasts.}
In practice return forecasts are always proprietary and come from sophisticated prediction techniques based on a wide variety of data available to the trading companies.
In this example, we simply add zero-mean noise to the realized returns to obtain the estimates and then rescale them to have a realistic mean square error of our prediction.
While these estimates are not real because they use the actual returns, they provide realistic values for the purpose of our computations.
We assume to know the risk-free interest rates exactly with $(\hat{r}_t)_{n+1} = (r_t)_{n+1}$.
The return estimates for the non-cash assets are $(\hat{r}_t)_{1:n} = \alpha((r_t)_{1:n} + \epsilon_t)$ where $\epsilon_t \sim \mathcal{N}(0, \sigma_\epsilon I)$ and $\alpha > 0$ is the scaling to minimize the mean-squared error $\Expect((\hat{r}_t)_{1:n} - (r_t)_{1:n})^2$~\cite[\Sec 7.3]{boyd2017}.
This method gives us return forecasts in the order of $\pm 0.3\%$ with an information ratio $\sqrt{\alpha} \approx 0.15$ typical of a proprietary return forecast prediction.

\paragraph{Parameters.}
The problem parameters are $\theta = (w_{t-1}, r_t, D_t, \Sigma^F_t, F_t)$ which, respectively, correspond to the previous assets allocation, vector of returns and the risk model matrices. Note that the risk model is updated at the beginning of every month.
Since the parameters not only affect the data in problem vectors, but also in the matrices, we cannot exploit the offline factorization caching for linear system~\eqref{eq:kkt_online}.

\paragraph{Problem setup.}
We simulate the trading system using the S\&P100 values~\cite{QuandlWIKI} from 2008 to 2013 with risk cost $\gamma = 100$, borrow cost $s_t = 0.0001$ and trade cost $\lambda = 0.01$. These values are similar as the benchmark values in~\cite{boyd2017}. We use different sparsity levels $c$ from $1$ to $40$.
Afterwards we collect data by sampling around the trajectory points from a hypersphere with radius $0.001$ times the magnitude of the parameter. For example, for a vector of returns of magnitude $\bar{r}_t$ we sample around $r_t$ with a radius $0.001\bar{r}_t$.
Even though this problem does not have to be solved in real-time, in order to optimize the trading performance, we must perform multiple expensive backtesting simulations.

\begin{figure}
	\begin{center}
	\small
\ref{perfportfolio}

\begin{tikzpicture}
	\begin{axis}
		[
		axis x line=bottom,
		axis y line=left,
		enlarge x limits=upper,  
		enlarge y limits=auto,
		ymode=log,          
		width=.45\textwidth,
		height=0.2\textheight,
		xlabel=$k$,
		ymin=1e-05,
		ymax=1e1,
		ylabel={Computation time [s]},
		xtick = data,
		cycle list name=black white,
		thick,
		legend columns=-1,
		legend to name=perfportfolio,
		legend entries = {$c=10$, $c=20$, $c=30$, $c=40$},
		]
		\foreach \t in {10, 20, 30, 40}{
			\addplot+
				table [x=n_best,
					y=mean_time_pred,
					col sep=comma] {./data/benchmarks/portfolio/portfolio_\t_performance.csv};
		}
	\end{axis}
\end{tikzpicture}
\begin{tikzpicture}
	\begin{axis}
		[
		axis x line=bottom,
		axis y line=left,
		enlarge x limits=upper,  
		enlarge y limits=auto,
		ymode=log,          
		width=.45\textwidth,
		height=0.2\textheight,
		xlabel=$k$,
		ylabel={Suboptimality},
		ymin=1e-05,
		ymax=1,
		xtick = data,
		cycle list name=black white,
		thick,
		]
		\foreach \t in {10, 20, 30, 40}{
			\addplot+ [error bars/.cd, y explicit, y dir=both]
				table [x=n_best,
					y=avg_subopt,
					col sep=comma] {./data/benchmarks/portfolio/portfolio_\t_performance.csv};
		}
	\end{axis}
\end{tikzpicture}

\begin{tikzpicture}
	\begin{axis}
		[
		axis x line=bottom,
		axis y line=left,
		enlarge x limits=upper,  
		enlarge y limits=auto,
		ymode=log,          
		width=.45\textwidth,
		height=0.2\textheight,
		xlabel=$k$,
		ylabel={Infeasibility},
		ymax=1,
		ymin=1e-06,
		xtick = data,
		cycle list name=black white,
		thick,
		]
		\foreach \t in {10, 20, 30, 40}{
			\addplot+
				table [x=n_best,
					y=avg_infeas,
					col sep=comma] {./data/benchmarks/portfolio/portfolio_\t_performance.csv};
		}
	\end{axis}
\end{tikzpicture}
\begin{tikzpicture}
	\begin{axis}
		[
		axis x line=bottom,
		axis y line=left,
		enlarge x limits=upper,  
		enlarge y limits=auto,
		width=.45\textwidth,
		height=0.2\textheight,
		ymax=110,
		ymin=0,
		xlabel=$k$,
		ylabel={Accuracy [\%]},
		xtick = data,
		cycle list name=black white,
		thick,
		]
		\foreach \t in {10, 20, 30, 40}{
			\addplot table [x=n_best, y=accuracy, col sep=comma] {./data/benchmarks/portfolio/portfolio_\t_performance.csv};
		}
	\end{axis}
\end{tikzpicture}
\end{center}
\caption{MLOPT average performance indicators for portfolio example.}
\label{fig:portfolio_mlopt}
\end{figure}
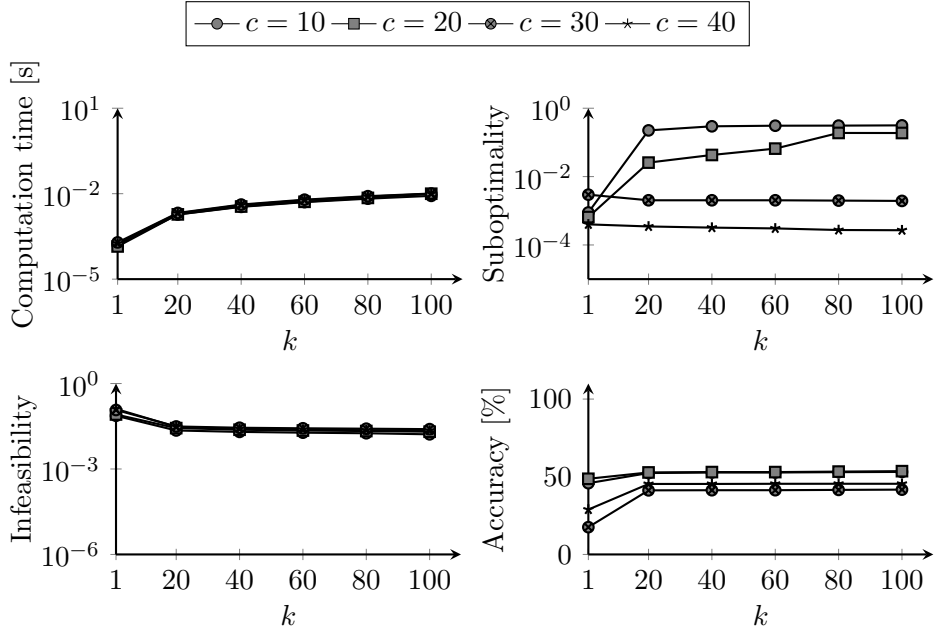

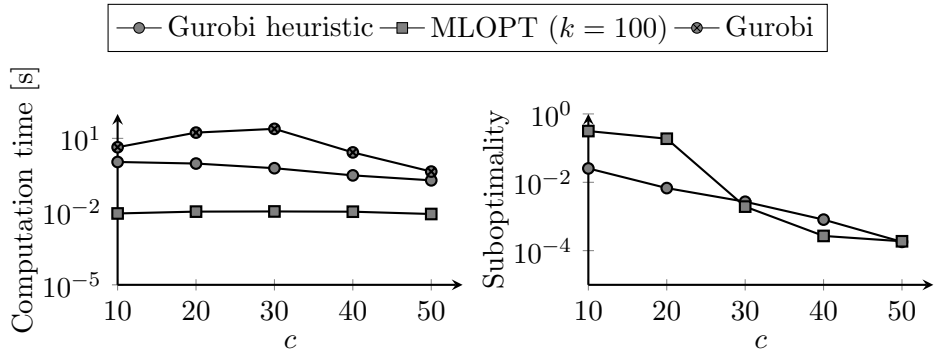
\begin{figure}
	\begin{center}
	\small
\ref{gurobiportfolio}

\begin{tikzpicture}
	\begin{axis}
		[
		axis x line=bottom,
		axis y line=left,
		enlarge x limits=upper,  
		enlarge y limits=auto,
		ymode=log,          
		width=.45\textwidth,
		height=0.2\textheight,
		xlabel=$c$,
		ylabel={Computation time [s]},
		ymin=1e-05,
		ymax=1e2,
		xtick = data,
		cycle list name=black white,
		thick,
		legend columns=-1,
		legend to name=gurobiportfolio,
		legend entries = {Gurobi heuristic, MLOPT $(k=100)$, Gurobi},
		]
		\addplot+ 
			table [x=k,
				y=mean_time_heuristic,
				col sep=comma] {./data/benchmarks/portfolio/portfolio_performance_full_heuristic.csv};
		\addplot+
			table [x=k,
				y=mean_time_pred,
				col sep=comma] {./data/benchmarks/portfolio/portfolio_performance_full_heuristic.csv};
		\addplot+ 
			table [x=k,
				y=mean_time_full,
				col sep=comma] {./data/benchmarks/portfolio/portfolio_performance_full_heuristic.csv};
	\end{axis}
\end{tikzpicture}
\begin{tikzpicture}
	\begin{axis}
		[
		axis x line=bottom,
		axis y line=left,
		enlarge x limits=upper,  
		enlarge y limits=auto,
		ymode=log,          
		width=.45\textwidth,
		height=0.2\textheight,
		xtick = data,
		xlabel=$c$,
		ylabel={Suboptimality},
		ymin=1e-05,
		ymax=1,
		cycle list name=black white,
		thick,
		]
		\addplot+
			table [x=k,
				y=avg_subopt_heuristic,
				col sep=comma] {./data/benchmarks/portfolio/portfolio_performance_full_heuristic.csv};
		\addplot+
			table [x=k,
				y=avg_subopt,
				col sep=comma] {./data/benchmarks/portfolio/portfolio_performance_full_heuristic.csv};
	\end{axis}
\end{tikzpicture}
\end{center}
\caption{Comparison between Gurobi and MLOPT performance for the portfolio example. Average computation time and suboptimality.}
\label{fig:portfolio_gurobi}
\end{figure}

\paragraph{Results.}
Table~\ref{tab:portfolio} shows the problem dimensions and the maximum computation time needed with each technique.
\reviewChanges{Here, the strategy pruning is not able to reduce $M$. This happens because, for most of the sample-strategy reassignments, the resulting problem becomes infeasible or suboptimal.}
Figure~\ref{fig:portfolio_mlopt} displays the performance of MLOPT for varying values of the top-$k$ strategies and in terms of computation time, suboptimality, infeasibility and accuracy.
For smaller $c$, although the total number of integer variable combinations is lower, the problem is harder to solve for every technique.
In this example, the performance does not significantly increase above $k=20$. We suspect that further accuracy improvements should happen with $k>100$.
In addition, suboptimality is very low for $c=30$ and $c=40$.
The performance comparison between the different methods appears in Figure~\ref{fig:portfolio_gurobi}.
MLOPT shows up to three orders of magnitude speedups over Gurobi and Gurobi heuristic despite moderate suboptimality and infeasibility values, mostly for $c=10$ and $c=20$.
With these computation time speedups, backtesting time can be significantly reduced and multiple parameter simulations can be executed to tune the problem parameters while evaluating the performance on historical data. This is crucial to obtain high quality portfolio trades.



%

\begin{table}
  \centering
	\footnotesize
  \caption{Portfolio problem dimensions and maximum times.}
  \label{tab:portfolio}
%

  \begin{tabular}
	  {@{}lllll
    S[round-precision = 4, table-number-alignment=left,scientific-notation=false]
    S[table-format=2.4, round-precision = 4, table-number-alignment=left,scientific-notation=false]
    S[round-precision = 4, table-number-alignment=left,scientific-notation=false]
    @{}}
    \toprule
    $c$ & $n_{\rm var}$ & $n_{\rm constr}$ & {\makecell[lb]{$M$\\ (unpruned)}} & $M$ & {\makecell[lb]{$t_{\rm max}$\\ MLOPT [s]}} & {\makecell[lb]{$t_{\rm max}$\\ Gurobi [s]}}& {\makecell[lb]{$t_{\rm max}$\\ Gurobi \\heuristic [s]}}\\
    \midrule
    \begin{filecontents*}{./data/benchmarks/portfolio/portfolio_performance_full_heuristic.csv}
,problem,learner,n_var,n_constr,n_theta,n_strategies,n_strategies_unpruned,n_infeas,avg_infeas,std_infeas,max_infeas,avg_subopt,std_subopt,max_subopt,avg_subopt_heuristic,std_subopt_heuristic,max_subopt_heuristic,mean_solve_time_pred,std_solve_time_pred,mean_pred_time_pred,std_pred_time_pred,mean_time_pred,std_time_pred,max_time_pred,mean_time_full,std_time_full,max_time_full,mean_time_heuristic,std_time_heuristic,max_time_heuristic,k
5,problem,pytorch,552,902,1581,4371,4371,4608,0.0164827648728283,0.0245773079303224,0.0999499754909286,0.3152794215278136,3.215288245612831,73.21499008167167,0.0254433979968767,0.6311837709671765,62.160702759035885,0.0083298444274402,0.0003864349087543,5.908313076924661e-05,1.3552527156068804e-20,0.0083889275582094,0.0003864349087543,0.0119058639077467,4.32871363207979,3.742602515894424,55.80382704734802,1.0775770395491069,0.1850378940192822,2.365391969680786,10
5,problem,pytorch,552,902,1581,2887,2887,7413,0.0203269189461898,0.0275918295209137,0.1557978507616361,0.1883793421811062,1.0566481983946294,12.558667285763216,0.00677186554798,0.0363247381129129,2.8621456396598752,0.0099260697180133,0.0002326673111561,4.698067696352581e-05,1.3552527156068804e-20,0.0099730503949768,0.0002326673111561,0.0136528137305523,17.196199869115162,77.87933413374692,2529.3926119804382,0.9285408677032,0.3624720591270267,2.024555206298828,20
5,problem,pytorch,552,902,1581,3371,3371,9097,0.0250591394651312,0.0257965375396026,0.1375670375031826,0.001942737730165,0.0109883663705326,0.1146459945749255,0.0026997932158876,0.0150211907895925,0.4951040292350811,0.010067963860596,0.0002818033701359,6.378556434305562e-05,2.710505431213761e-20,0.0101317494249391,0.0002818033701359,0.0136586513632932,24.49181538257959,98.25585680733634,2313.78190779686,0.5987072604415197,0.4968593187228769,2.322695016860962,30
5,problem,pytorch,552,902,1581,3736,3736,9881,0.0222040257007261,0.0251318561112188,0.1966263921249615,0.0002696733168635,0.000781820102738,0.0037682576446571,0.0008094358242416,0.0057698674440218,0.1605304520093073,0.0098506576846822,0.0002092909194392,5.6273579715853794e-05,2.710505431213761e-20,0.0099069312643981,0.0002092909194392,0.0133702822924844,2.6555613181323494,9.553343101639866,226.24012899398804,0.3042984230160832,0.3917545191191642,1.6632061004638672,40
5,problem,pytorch,552,902,1581,4130,4130,9854,0.0203977252152025,0.0220604800221191,0.192783136421356,0.0001881034583674,0.000284560946832,0.0013354111925988,0.0001807153054205,0.0020736033525286,0.0778132949609061,0.0079168721249229,0.0002060599833443,5.05821470470864e-05,1.3552527156068804e-20,0.0079674542719699,0.0002060599833443,0.0113470928434582,0.4359126856689301,1.0226375184294343,19.668039083480835,0.1906458550673847,0.2990350040504249,1.4714548587799072,50
\end{filecontents*}
\csvreader[head to column names, late after line=\\]{./data/benchmarks/portfolio/portfolio_performance_full_heuristic.csv}{
    k=\c,
    n_constr=\nconstr,
    n_var=\nvar,
    n_strategies=\M,
    n_strategies_unpruned=\Munpruned,
    max_time_pred=\tmaxmlopt,
    max_time_full=\tmaxgurobi,
    max_time_heuristic=\tmaxheuristic,
    }{\c & \nvar & \nconstr & \Munpruned & \M & \tmaxmlopt & \tmaxgurobi & \tmaxheuristic}
\bottomrule
  \end{tabular}
\end{table}

\subsection{Motion Planning}%
\label{sub:motion}
We consider the problem of motion planning in the presence of obstacles.
This problem has a wide variety of applications including autonomous vehicles, space robots, and~\glspl{UAV}~\cite{planningalgorithms}.
These problems must usually be solved online within a few milliseconds to provide inputs that are frequent enough to control the system dynamics.
Unfortunately, the presence of obstacles makes the problem nonconvex and, therefore, very challenging to solve online.
In the literature, several approaches have been proposed to model the motion planning problem with~\gls{MIO}~\cite{milp_path}.
However, solution times in the order of a few milliseconds are still out of reach of state-of-the-art \gls{MIO} solvers.

We consider the problem in $d$ dimensions, in practice $d=2$ for planar systems or $d=3$ for aerial systems.
The state of the system is $x_t = (p_t, v_t)$ where $p_t\in \reals^d$ is the position at time $t$ and $v_t \in \reals^d$ the velocity.  The input $u_t \in \reals^d$ are the forces produced by the system in every direction.
For example, for an~\gls{UAV}, $u_t$ corresponds to the thruster forces.
The discrete-time linear system dynamics are described as
\begin{equation*}
	(p_{t+1}, v_{t+1}) = A(p_t, v_t) + Bu_t,\qquad t=0,\dots,T,
\end{equation*}
where  $A\in \reals^{2d \times 2d}$ and $B\in \reals^{2d \times d}$ and $\tau > 0$ is the sampling time.
The initial state of the system is $(p^{\rm init}, v^{\rm init})$.
We define upper and lower bounds on the state and inputs as
\begin{equation*}
	\begin{aligned}
	    &\underline{p} \le p_t \le \overline{p},\qquad \underline{v} \le v_t \le \overline{v},\qquad &t=0,\dots,T,\\
	    &\underline{u} \le u_t \le \overline{u},\qquad &t=0,\dots,T-1.
	\end{aligned}
\end{equation*}
We model every obstacle $i$ as a rectangle in $\reals^d$ with upper bounds $\overline{o}^i \in \reals^{d}$ and lower bounds $\underline{o}^i\in \reals^{d}$ for $i=1,\dots,n_{\rm obs}$.
At every time $t$ and for every obstacle $i$, we model the obstacle avoidance decisions with binary variables $\overline{\delta}_{t}^{i} \in \{0, 1\}^d$ for the obstacle upper bounds and $\underline{\delta}_{t}^{i} \in \{0, 1\}^d$ for the lower bounds.
We can write obstacle avoidance as the following big-M conditions
\begin{equation*}
	\overline{o}^i - M \overline{\delta}_{t}^{i} \le p_t \le \underline{o}^i + M \underline{\delta}_{t}^{i},\qquad t=0,\dots,T,\quad i=1,\dots,n_{\rm obs}.
\end{equation*}
In addition, we must impose that we cannot be at the same time on different sides of an obstacle,
\begin{equation*}
	\ones^\tpose \underline{\delta}_{t}^{i} + \ones^\tpose \overline{\delta}_{t}^{i}\le 2d - 1,\qquad t=0,\dots,T,\quad i=1,\dots,n_{\rm obs},
\end{equation*}
where $\ones$ is the vector of ones of dimension $d$.
Our goal is to minimize the distance with respect to a desired position $(p_t \approx p^{\rm des})$ while expending minimal control effort $(u_t \approx 0)$. We can write the cost as
\begin{equation*}
	\|p_T - p^{\rm des}\|_2^2 + \sum_{t=0}^{T-1} \|p_t - p^{\rm des}\|_2^2 + \gamma \|u_t\|^2_2,
\end{equation*}
where $\gamma > 0$ balances the trade-off between performance and control effort.
The motion planning problem can be written as
\begin{equation}
	\label{eq:trajectory}
	\begin{array}{ll}
		\text{minimize} & \displaystyle \|p_T - p^{\rm des}\|_2^2 + \sum_{t=0}^{T-1} \|p_t - p^{\rm des}\|_2^2 + \gamma \|u_t\|^2_2\\
		\text{subject to} & (p_{t+1}, v_{t+1}) = A(p_t, v_t) + Bu_t,\qquad t=0,\dots,T-1\\
				  & p_0 = p^{\rm init},\quad v_0 = v^{\rm init}\\
				  & \overline{o}^i - M \overline{\delta}_{t}^{i} \le p_t \le \underline{o}^i + M \underline{\delta}_{t}^{i},\qquad t=0,\dots,T,\quad i=1,\dots,n_{\rm obs}\\
				  & \ones^\tpose \underline{\delta}_{t}^{i} + \ones^\tpose \overline{\delta}_{t}^{i} \le 2d - 1\qquad t=0,\dots,T\\
				  & \overline{\delta}_t^i \in \{0, 1\}^d,\quad \underline{\delta}_t^i \in \{0, 1\}^d,\qquad i=1,\dots,n_{\rm obs}\\
				  & \reviewChanges{p_t \in [\underline{p}, \overline{p}],\; v_t \in [\underline{v}, \overline{v}], \qquad t=0,\dots,T,}\\
				  & u_t \in [\underline{u}, \overline{u}],\qquad t=0,\dots,T-1.\\
	\end{array}
\end{equation}

\paragraph{Parameters.}
We consider the problem of computing the optimal trajectory from any point in the plane to the desired position $p^{\rm des}$.
In this example the rest of the problem data does not change. Therefore, the problem parameters are the initial state $\theta = p^{\rm init}$.
Problem~\eqref{eq:trajectory} must be solved within each sampling time $\tau$ to provide fast enough inputs to the dynamical system.

\begin{figure}[h]
	\centering
	\includegraphics[width=0.7\linewidth]{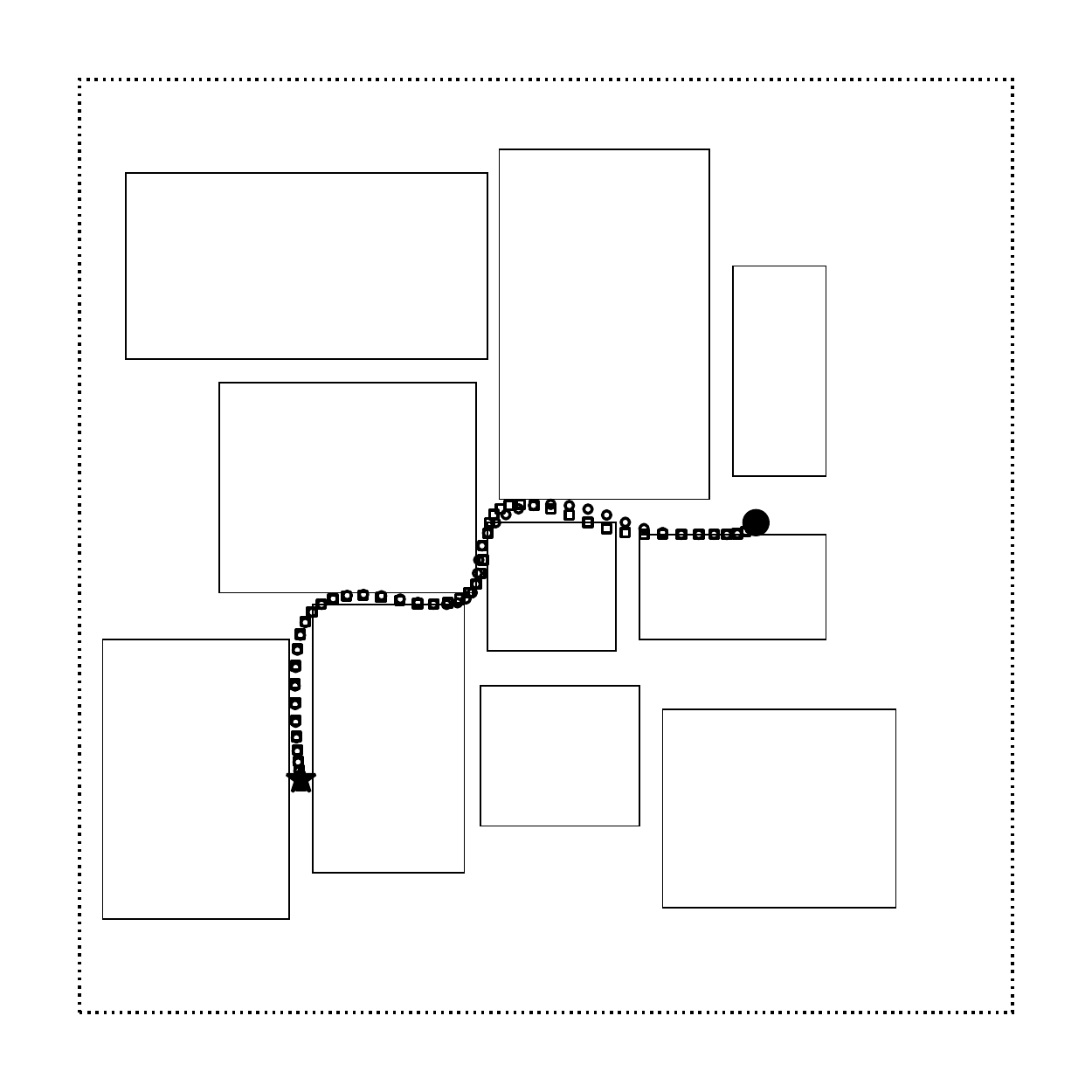}
	\caption{Trajectory planning example with 10 obstacles and sampling time 0.1 sec. Circles indicate the optimal path computed with Gurobi in 18.18 sec and the squares indicate the optimal path computed with MLOPT with $k=100$ in 0.06 sec. The filled dot is the starting position while the \reviewChanges{start} is the desired position. The boundaries of the feasible positions are the dotted line while obstaces are displayed as rectangles.}%
	\label{fig:obstacle}
\end{figure}

\paragraph{Problem setup.}
We consider a discrete-time  double-integrator dynamical system in two dimensions as in~\cite{milp_path}.
The optimizer must provide solutions within the sampling time $\tau =0.1\; {\rm sec}$ in order to operate the system properly.
The time horizon is $T=60$ which corresponds to $6\; {\rm sec}$.
The cost tradeoff parameter is $\gamma=0.01$.
The initial state is given by $p^{\rm init} = \theta$ and $v^{\rm init } = 0$.
The desired position is $p^{\rm des} = (-10.5, -10)$ and the desired velocity $v^{\rm des} = 0$.
Figure~\ref{fig:obstacle} shows an example trajectory for $\theta = (9, 13)$.
To generate training samples $\theta_i$, we sample uniformly between $\underline{p}$ and $\overline{p}$.

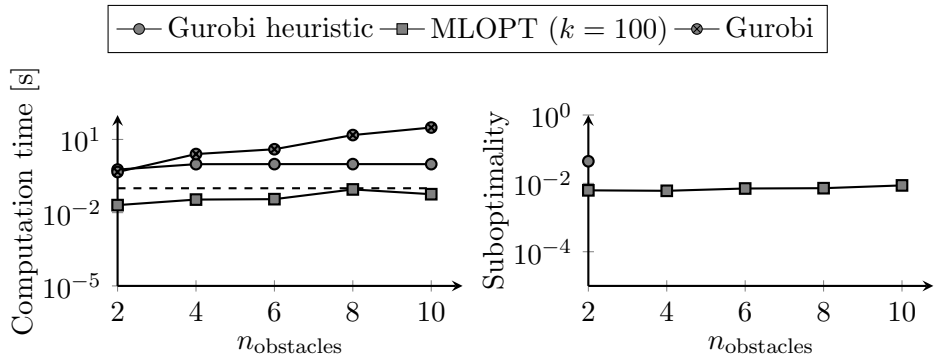
\begin{figure}
	\begin{center}
	\small
\ref{gurobiobstacle}

\begin{tikzpicture}
	\begin{axis}
		[
		axis x line=bottom,
		axis y line=left,
		enlarge x limits=upper,  
		enlarge y limits=auto,
		ymode=log,          
		width=.45\textwidth,
		height=0.2\textheight,
		xlabel=$n_{\rm obstacles}$,
		ylabel={Computation time [s]},
		ymin=1e-05,
		ymax=1e2,
		xtick = data,
		cycle list name=black white,
		thick,
		legend columns=-1,
		legend to name=gurobiobstacle,
		legend entries = {Gurobi heuristic, MLOPT $(k=100)$, Gurobi},
		]
		\addplot+
			table [x=n_obs,
				y=mean_time_heuristic,
				col sep=comma] {./data/benchmarks/obstacle/obstacle_performance_full_heuristic.csv};
		\addplot+ 
			table [x=n_obs,
				y=mean_time_pred,
				col sep=comma] {./data/benchmarks/obstacle/obstacle_performance_full_heuristic.csv};
		\addplot+ 
			table [x=n_obs,
				y=mean_time_full,
				col sep=comma] {./data/benchmarks/obstacle/obstacle_performance_full_heuristic.csv};
     \addplot[mark=none, dashed, thick] coordinates {(2,0.1) (10,0.1)};  
	\end{axis}
\end{tikzpicture}
\begin{tikzpicture}
	\begin{axis}
		[
		axis x line=bottom,
		axis y line=left,
		enlarge x limits=upper,  
		enlarge y limits=auto,
		ymode=log,          
		width=.45\textwidth,
		height=0.2\textheight,
		xtick = {2, 4, 6, 8, 10},
		xlabel=$n_{\rm obstacles}$,
		ylabel={Suboptimality},
		ymin=1e-05,
		ymax=1,
		cycle list name=black white,
		thick,
		]
		\addplot+
			table [x=n_obs,
				y=avg_subopt_heuristic,
				col sep=comma] {./data/benchmarks/obstacle/obstacle_performance_full_heuristic.csv};
		\addplot+
			table [x=n_obs,
				y=avg_subopt,
				col sep=comma] {./data/benchmarks/obstacle/obstacle_performance_full_heuristic.csv};
	\end{axis}
\end{tikzpicture}
\end{center}
\caption{Comparison between Gurobi and MLOPT performance for the motion planning example. Average computation time and suboptimality. The dashed line indicates the sampling time.\reviewChanges{Missing points indicate problems for which Gurobi heuristic is not able to return a feasible solution within the maximum time $1$ sec.}}
\label{fig:obstacle_gurobi}
\end{figure}

\paragraph{Results.}
Table~\ref{tab:obstacle} shows the problem dimensions along with the maximum computation time taken by each technique.
\reviewChanges{Here, the strategy pruning is not able to significantly reduce $M$. This happens because, in this case, the samples already provide an already small, although representative, selection of unique strategies to use.}
Figure~\ref{fig:portfolio_mlopt} shows the average performance of MLOPT for varying top-$k$ strategies selected.
Even if the accuracy is always low, most of the times less than 10\%, the suboptimality and infeasibility are acceptable for this application.
As shown in Figure~\ref{fig:obstacle}, there can be some constraints violations due to marginal violations of the obstacle boundaries.
We show the direct comparison of Gurobi, Gurobi heuristic and MLOPT in Figure~\ref{fig:obstacle_gurobi}. The time improvements of MLOPT allow its real-time implementation, in contrast to state-of-the-art options.
\reviewChanges{Remarkably, except in the smallest case of $n_{\rm obs}=2$, Gurobi heuristic is not able to return a solution within the maximum time we allowed, i.e., $1$ sec. Note that this limit is already $10$ times higher than the sampling time $\tau = 0.1$ sec. }

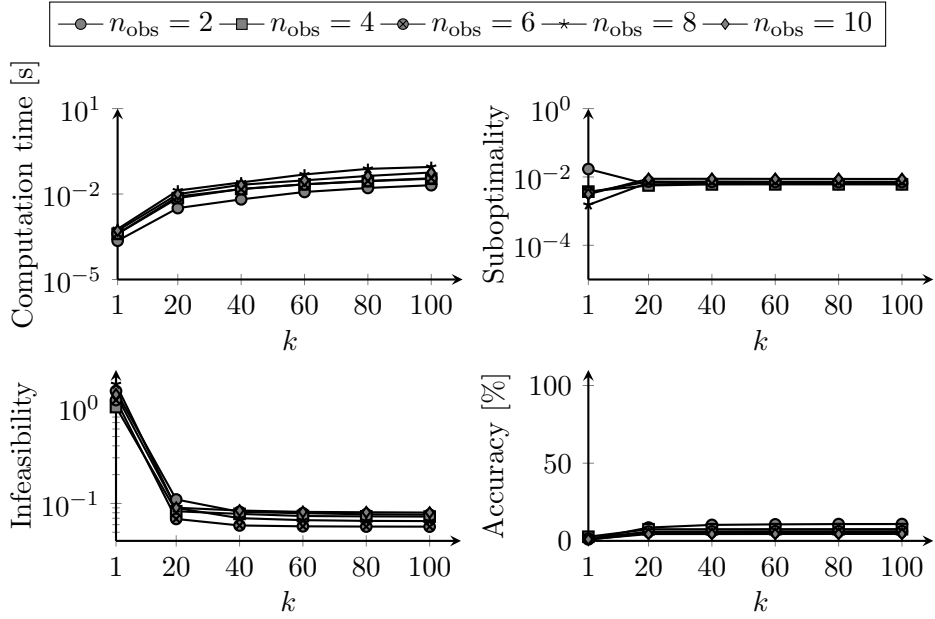
\begin{figure}
	\begin{center}
	\small
\ref{perfobstacle}

\begin{tikzpicture}
	\begin{axis}
		[
		axis x line=bottom,
		axis y line=left,
		enlarge x limits=upper,  
		enlarge y limits=auto,
		ymode=log,          
		width=.45\textwidth,
		height=0.2\textheight,
		xlabel=$k$,
		ymin=1e-05,
		ymax=1e1,
		ylabel={Computation time [s]},
		xtick = data,
		cycle list name=black white,
		thick,
		legend columns=-1,
		legend to name=perfobstacle,
		legend entries = {$n_{\rm obs}=2$, $n_{\rm obs}=4$, $n_{\rm obs}=6$, $n_{\rm obs}=8$,$n_{\rm obs}=10$}
		]
		\foreach \t in {2,4,6,8,10}{
			\addplot+
				table [x=n_best,
					y=mean_time_pred,
					col sep=comma] {./data/benchmarks/obstacle/obstacle_\t_performance.csv};
		}
	\end{axis}
\end{tikzpicture}
\begin{tikzpicture}
	\begin{axis}
		[
		axis x line=bottom,
		axis y line=left,
		enlarge x limits=upper,  
		enlarge y limits=auto,
		ymode=log,          
		width=.45\textwidth,
		height=0.2\textheight,
		xlabel=$k$,
		ymin=1e-05,
		ymax=1,
		ylabel={Suboptimality},
		xtick = data,
		cycle list name=black white,
		thick,
		]
		\foreach \t in {2,4,6,8,10}{
			\addplot+ 
				table [x=n_best,
					y=avg_subopt,
					col sep=comma] {./data/benchmarks/obstacle/obstacle_\t_performance.csv};
		}
	\end{axis}
\end{tikzpicture}

\begin{tikzpicture}
	\begin{axis}
		[
		axis x line=bottom,
		axis y line=left,
		enlarge x limits=upper,  
		enlarge y limits=auto,
		ymode=log,          
		width=.45\textwidth,
		height=0.2\textheight,
		xlabel=$k$,
		ylabel={Infeasibility},
		xtick = data,
		cycle list name=black white,
		thick,
		]
		\foreach \t in {2,4,6,8,10}{
			\addplot+
				table [x=n_best,
					y=avg_infeas,
					col sep=comma] {./data/benchmarks/obstacle/obstacle_\t_performance.csv};
		}
	\end{axis}
\end{tikzpicture}
\begin{tikzpicture}
	\begin{axis}
		[
		axis x line=bottom,
		axis y line=left,
		enlarge x limits=upper,  
		enlarge y limits=auto,
		width=.45\textwidth,
		height=0.2\textheight,
		ymax=110,
		ymin=0,
		xlabel=$k$,
		ylabel={Accuracy [\%]},
		xtick = data,
		cycle list name=black white,
		thick,
		]
		\foreach \t in {2,4,6,8,10}{
			\addplot table [x=n_best, y=accuracy, col sep=comma] {./data/benchmarks/obstacle/obstacle_\t_performance.csv};
	        }
	\end{axis}
\end{tikzpicture}
\end{center}
\caption{MLOPT average performance indicators for motion planning example. The dashed line indicates the sampling time.}
\label{fig:obstacle_mlopt}
\end{figure}

\begin{table}
	\footnotesize
  \centering
  \caption{Motion planning problem dimensions and maximum times.}
  \label{tab:obstacle}
  \begin{tabular}
	  {@{}lllll
    S[round-precision = 4, table-number-alignment=left,scientific-notation=false]
    S[table-format=2.4, round-precision = 4, table-number-alignment=left,scientific-notation=false]
    S[round-precision = 4, table-number-alignment=left,scientific-notation=false]
    @{}}
    \toprule
    $n_{\rm obs}$ & $n_{\rm var}$ & $n_{\rm constr}$ & {\makecell[lb]{$M$\\ (unpruned)}} & $M$ & {\makecell[lb]{$t_{\rm max}$\\ MLOPT [s]}} & {\makecell[lb]{$t_{\rm max}$\\ Gurobi [s]}}& {\makecell[lb]{$t_{\rm max}$\\ Gurobi \\heuristic [s]}}\\
    \midrule
    \begin{filecontents*}{./data/benchmarks/obstacle/obstacle_performance_full_heuristic.csv}
,problem,learner,n_var,n_constr,n_theta,n_strategies,n_strategies_unpruned,n_infeas,avg_infeas,std_infeas,max_infeas,avg_subopt,std_subopt,max_subopt,avg_subopt_heuristic,std_subopt_heuristic,max_subopt_heuristic,mean_solve_time_pred,std_solve_time_pred,mean_pred_time_pred,std_pred_time_pred,mean_time_pred,std_time_pred,max_time_pred,mean_time_full,std_time_full,max_time_full,mean_time_heuristic,std_time_heuristic,max_time_heuristic,n_obs
5,problem,pytorch,1135,3773,2,1351,1357,8533,0.0763337241814911,0.1198804967063059,2.070119922567559,0.0062525716258529,0.0264371069816416,0.2203458151170941,0.0442972123260463,0.1882809259841517,5.820870983192574,0.0203785367727279,0.0011007351660845,2.5644302368164053e-05,1.0164395367051604e-20,0.0204041810750961,0.0011007351660845,0.0262662315368652,0.4621566620111465,0.4718905854328575,3.531091928482056,0.5818965394973755,0.4335641022379973,1.0085110664367676,2
5,problem,pytorch,1615,10133,2,1112,1114,8939,0.0729244817596623,0.1367267874411599,1.655217710736871,0.0060619980027848,0.0194252420697003,0.1838672199231256,inf,,inf,0.0343296449661254,0.0040786652425041,4.9859976768493655e-05,0.0,0.0343795049428939,0.0040786652425041,0.075171262717247,2.517781358075142,1.891744176213433,12.793086051940918,0.9685106461763382,0.1538714877781156,1.0359549522399902,4
5,problem,pytorch,2095,20333,2,917,921,9145,0.057234594398359,0.082736284565508,1.3783633985879338,0.0070785956626373,0.0213600594199911,0.1838672199231256,inf,,inf,0.035737982392311,0.0052097495528204,1.5633463859558106e-05,0.0,0.0357536158561706,0.0052097495528204,0.0503505638837814,3.960185301613808,3.7575087775248712,248.4688320159912,0.9741608147859572,0.1568967540147091,1.0308640003204346,6
5,problem,pytorch,2575,34373,2,818,818,9302,0.0654274702655966,0.1097134848827221,1.5884871122083977,0.0072295356195482,0.0209504777052356,0.1838672199231254,inf,,inf,0.0894406996488571,0.0191543511996907,1.6277360916137695e-05,0.0,0.0894569770097732,0.0191543511996907,0.2213493042469024,15.071509842181204,22.687043278172517,395.543251991272,0.9745815408468248,0.1540735725867564,1.036055088043213,8
5,problem,pytorch,3055,52253,2,657,658,9344,0.0802849039877339,0.1337975310718714,1.48440411536463,0.0087183769273381,0.0219501061720454,0.179928839122721,inf,,inf,0.0570935250759124,0.0025260387592714,1.1868309974670412e-05,1.694065894508601e-21,0.0571053933858871,0.0025260387592714,0.0613791416406631,30.728369207954408,47.861150042147344,809.3376410007477,0.971151594543457,0.1566677397870518,1.0465240478515625,10
\end{filecontents*}
\csvreader[head to column names, late after line=\\]{./data/benchmarks/obstacle/obstacle_performance_full_heuristic.csv}{
    n_obs=\nobs,
    n_constr=\nconstr,
    n_var=\nvar,
    n_strategies=\M,
    n_strategies_unpruned=\Munpruned,
    max_time_pred=\tmaxmlopt,
    max_time_full=\tmaxgurobi,
    max_time_heuristic=\tmaxheuristic,
    }{\nobs & \nvar & \nconstr & \Munpruned & \M & \tmaxmlopt & \tmaxgurobi & \tmaxheuristic}
\bottomrule
  \end{tabular}
\end{table}

\subsection{Remarks}%
\label{sub:remarks}
The numerical examples show several benefits of our approach.
First, the strategy pruning technique allows us to work with a total number of strategies between 100 and 10,000.
These numbers allow high quality predictions in terms of suboptimality and infeasibility with \glspl{NN} combined with top-$k$ strategies selection.
Second, despite low accuracy when comparing the predicted strategy to the exact optimal one, the average \reviewChanges{suboptimality} is comparable or better than the one of Gurobi heuristic.
In addition, infeasibility is always within acceptable levels for the proposed application.
This means that, even in case of multiple global minimizers, our method is able to consistently predict optimal or close-to-optimal solutions.
Third, \reviewChanges{for the problems considered,} time measurements show up to three orders of magnitude speedups of MLOPT compared to Gurobi and Gurobi heuristics, both in terms of average and worst-case solution time.
Since for the fuel cell energy management and the motion planning examples, optimal solutions must be computed within hard real-time requirements, only fast optimization methods can be implemented in practice and Gurobi-based approaches are too slow. On the contrary, MLOPT proves to be a much faster approach for computing online solutions.

\section{Conclusions}
\label{sec:conclusions}

We proposed a machine learning method for solving online~\gls{MIO} at very high speed.
By using the Voice of Optimization framework~\cite{bertsimas2018} we exploited the repetitive nature of online optimization problems to greatly reduce the solution time.
Our method casts the core part of the optimization algorithm as a multiclass classification problem that we solve using a \gls{NN}.
In this work we considered the class of~\gls{MIQO} which, while covering the vast majority of real-world optimization problems, allows us to further reduce the computation time exploiting its structure.
For \gls{MIQO}, we only have to solve a linear system to obtain the optimal solution after the \gls{NN} prediction.
In other words, our approach does not require any solver nor iterative routine to solve parametric \glspl{MIQO}.
Our method is not only extremely fast but also reliable. Compared to branch-and-bound methods  our approach has a very predictable online solution time for which we can exactly bound the~\glspl{flop} complexity.
This is of crucial importance for real-time applications.
We benchmarked our approach against the state-of-the-art solver Gurobi on different real-world examples showing $100$ to $1000$ fold speedups in online execution time.

There are several future research directions to improve the MLOPT framework.
First, specialized neural networks architectures might provide more accurate strategy classifiers for certain classes of optimization problems.
Second, the MLOPT training phase requires computing the optimal solution of problems with several different parameter combinations. Therefore, the current approach does not apply to problems that take hours to solve with Gurobi.
Using suboptimal solutions in the training phase would allow us significantly reduce the training time and apply MLOPT to such problems.
Third, strategy reduction techniques combined with structure exploiting classifiers would enable applications to large-scale optimization problems with only discrete variables.
Finally, including the knowledge about the previous solution would improve the performance of MLOPT, similarly to how warm-starting techniques help reducing computations in~\gls{MIO} solvers.

\ifpreprint
\section*{Acknowledgments}
\reviewChanges{
The authors would like to thank the anonymous reviewers for providing comments that greatly helped improving the quality of this paper.

The simulations presented in this article were performed on computational resources managed and supported by Princeton Research Computing, a consortium of groups including the Princeton Institute for Computational Science and Engineering (PICSciE) and the Office of Information Technology's High Performance Computing Center and Visualization Laboratory at Princeton University.}

\fi


\bibliography{refs}

\end{document}